\documentclass[reqno,11pt]{amsart}
\usepackage{amsfonts, amsmath, amssymb, amsthm, calrsfs, color}
  \usepackage{amsfonts}
\usepackage{amsmath}
\allowdisplaybreaks[1]

\def\bbm[#1]{\mbox{\boldmath $#1$}}

\newtheorem{thm}{Theorem}[section]
\newtheorem{lemma}[thm]{Lemma}
\newtheorem{prop}[thm]{Proposition}
\newtheorem{cor}[thm]{Corollary}

\theoremstyle{definition}
 \numberwithin{equation}{section}

   \newcommand{\C}{\mathbb{C}}

   \newcommand{\N}{\mathbb{N}}
\newcommand{\D}{\Delta}
\newcommand{\rr}{\mathbb{R}}
\newcommand{\intr}{\int_{\R^2}}
\newcommand{\R}{\mathbb{R}}
\newcommand{\al}{\alpha}
\newcommand{\de}{\delta}
 
\newcommand{\la}{\lambda}
\newcommand{\into}{\int_\Omega}

 \newcommand{\e}{\varepsilon}
 \renewcommand{\(}{\left(}
\renewcommand{\)}{\right)}

\newcommand{\beq}{\begin{equation}}
\newcommand{\eeq}{\end{equation}}
\def\bbm[#1]{\mbox{\boldmath $#1$}}
\def\Re{{\rm Re}}
\def\Im{{\rm Im}}

\begin{document}

\title[Liouvillle equation with
a singular source]{Blow-up phenomena for the Liouville equation \\with
a singular source of integer multiplicity}
\author{Teresa D'Aprile}
\address[Teresa D'Aprile] {Dipartimento di Matematica, Universit\`a di Roma ``Tor
Vergata", via della Ricerca Scientifica 1, 00133 Roma, Italy.}
\email{daprile@mat.uniroma2.it}



\begin{abstract}
We are concerned with the existence of blowing-up solutions 
 to the following boundary value problem $$-\Delta u= \la a(x) e^u-4\pi N \delta_0\;\hbox{ in } \Omega,\quad	u=0 \;\hbox{ on }\partial \Omega,$$  where $\Omega$ is a smooth  and bounded domain in $\R^2$ such that $0\in\Omega$,  $a(x)$ is a positive smooth function, $N$ is a  positive integer  and $\la>0$ is a small parameter. Here $\delta_0$ defines the Dirac measure with pole at $0$.   We find conditions on the function $a$ and on the domain $\Omega$ under which  
there exists a solution $u_\la$ blowing up at $0$ and 
satisfying $\la\into a(x)e^{u_\la} \to 8\pi(N+1)$ as $\la\to 0^+$. 

\medskip

\noindent {\bf Mathematics Subject Classification 2010:} 35J20, 35J57,
35J61

\noindent {\bf Keywords:} singular Liouville equation, blowing-up solutions,
perturbation methods

\end{abstract}

\maketitle

\section{Introduction}

Let $\Omega$ be a bounded domain in $\R^2$ with a smooth boundary containing the origin. In this paper we consider the following Liouville equation with Dirac mass measure  


\begin{equation}\label{eq}
  \left\{
      \begin{aligned}&- \D u = \la a(x) e^u-4\pi N \delta_0&  \hbox{ in }&  \Omega,\\
    &  \ u=0 &  \hbox{ on }& \partial \Omega.
  \end{aligned}
    \right. \end{equation}
Here $\la$ is a positive small parameter,  $\delta_0$ denotes Dirac mass supported at $0$, $a$ is a  smooth function satisfying $\inf_\Omega a(x)>0$  and $N $ is a positive integer.

Problem \eqref{eq} is motivated by its links with the modeling of physical phenomena. In particular,  \eqref{eq} arises in the study of vortices   in a planar model of Euler flows (see  \cite{delespomu}, \cite{bapi}). In vortex theory the interest in constructing \textit{blowing-up} solutions is related to relevant physical properties, in particular  the presence of vortices with a strongly localised electromagnetic field.

The asymptotic behaviour of  family 
of blowing up solutions 
     can be referred to the papers \cite{bapa}, \cite{breme}, \cite{lisha}, \cite{mawe},  \cite{nasu}, \cite{su} for the regular problem, i.e. when $N=0$. An extension to the singular case $N>0$
      is contained in \cite{bachenlita}-\cite{bata}. 

 The   analysis  of the blowing-up behaviour at points away from $0$ actually is very similar to the asymptotic  analysis arising in the regular case which  has been pursued with success   and, at the present time, an accurate description of the  concentration  phenomenon is available. Precisely, the 
 analysis in the above  works yields that if $u_\la$ is    an unbounded family of solutions  of \eqref{eq}   for which   $\la\into a(x) e^{u_\la} $ is uniformly bounded and $u_\la$ is uniformly bounded in a neighborhood of $0$, then, up to a subsequence,  there is an integer $m\geq 1$ such that 
   \beq\label{issue}\la\into a(x)e^{u_\la}dx\to 8\pi m\hbox{ as }\la\to 0^+.\eeq Moreover there are points $\xi_1^\la,\ldots, \xi_m^\la\in\Omega$ which remain uniformly distant from the boundary $\partial\Omega$, from $0$  and from one another such that  \beq\label{issue11}\la a(x)e^{u_\la}-8\pi\sum_{j=1}^m \delta_{\xi_j^\la}\to 0\eeq in the measure sense. 
  Also  the location of the     blowing-up points is well understood when concentration occurs away from $0$.
    Indeed, in \cite{nasu} and \cite{su} it is established that the $m$-tuple $(\xi_1^\la,\ldots, \xi_m^\la)$ converges, up to a subsequence, to a critical point of the functional 
  \beq\label{natur}\frac12\sum_{j=1}^mH(\xi_j,\xi_j)+\frac12\sum_{j,h=1\atop j\neq h}^mG(\xi_j,\xi_h)-\frac{N}{2}\sum_{j=1}^mG(\xi_j,0).\eeq  
 Here 
   $G(x,y)$ is the Green's function of $-\Delta$ over $\Omega$ under Dirichlet
boundary conditions and  $H(x,y)$ denotes its regular part:
$$H(x,y)=G(x,y)-\frac{1}{2\pi}\log\frac{1}{|x-y|}.$$
The above description of blowing-up behaviour  continues to work if we are  in the presence of multiples singularities $\sum_i N_i\delta_{p_i}$ in \eqref{eq}, provided that we substitute the term  $  \frac{N}{2}\sum_jG(\xi_j,0)$ by $\sum_{i}\frac{N_i}{2}\sum_jG(\xi_j,p_i)$ in \eqref{natur}.
\ 

  The reciprocal issue, namely the existence of positive solutions with the property  \eqref{issue11}, has been addressed for the regular case $N=0$ first   in \cite{we} in the case of  a single point of concentration (i.e. $m=1$), later generalised   to the case of multiple concentration associated to any nondegenerate critical point of the functional \eqref{natur} (\cite{bapa}, \cite{chenlin}) or, more generally, to a any \textit{topologically nontrivial} critical point (\cite{delkomu}-\cite{espogropi}). 
 In particular, still for $N=0$, a family of solutions $u_\la$ concentrating at  $m$-tuple of points as $\la\to 0^+$  has been found in some special cases: for any $m\geq 1$, provided that  $\Omega$ is not simply connected  (\cite{delkomu}), and  for $m\in\{1,\ldots, h\}$ if $\Omega$ is a $h$-dumbell with thin handles (\cite{espogropi}). We mention that  functionals similar to  \eqref{natur}  occur to detect multiple-bubbling solutions in different contexts, see \cite{bapiwe}, \cite{espomupi1}, \cite{espomupi2}, \cite{weiyezhou} for other related singularly perturbed problems.

   In the singular case $N>0$    solutions which concentrate in the measure sense at  $m$ distinct points away from $0$ have been built in \cite{delkomu} provided that $m<1+N$. This result has been extended in \cite{daprile}  to
the case of multiple singular sources: in particular it is showed  that, under suitable
restrictions on the weights, if several sources exist then the more involved topology
generates a large number of blow-up solutions.

We point out that in all the above results  
   concentration occurs 
    at points different from the location of the source. The problem of finding solutions with additional concentration around the source is of different nature.   In case they exist, the blowing-up at the singularity  provides an additional contribution of $8\pi(1+N)\delta_0$ in the limit  \eqref{issue}, see \cite{bachenlita},  \cite{bata}, \cite{espo}, \cite{ta1}, \cite{ta2}. 
More precisely the asymptotic analysis in the general case can be formulated as follows:  if $u_\la$ is    an unbounded family of solutions  of \eqref{eq}   for which   $\la\into a(x)e^{u_\la} $ is uniformly bounded and $u_\la$ is unbounded in any neighborhood of $0$, then, up to a subsequence,  there is an integer $m\geq 0$ such that 
   $$\la\into a(x)e^{u_\la}dx\to 8\pi m+8\pi(N+1)\hbox{ as }\la\to 0^+.$$ Moreover there are $m$ distinct points $\xi_1,\ldots, \xi_m\in\Omega\setminus\{0\}$  such that, up to  subsequence,    \beq\label{issue1}\la a(x)e^{u_\la}\to 8\pi\sum_{j=1}^m \delta_{\xi_j}+8\pi(N+1)\delta_0\eeq in the measure sense.   
We mention that also in this case the analysis can be generalized to any number of sources. Moreover, under some extra assumptions it is possible to define a functional which replaces  \eqref{natur} in locating the points $\xi_j$ where the concentration occurs, anyway to avoid technicalities we  will not go into any further detail (see \cite{espo}). 

  The question on the existence of  solutions     to \eqref{eq} concentrating at $0$     is far from being completely settled. Indeed only partial results are known: in \cite{espo} the construction of solutions concentrating at $0$  is carried out provided that  $N\in(0,+\infty)\setminus \N$. 
 To our knowledge, the only paper dealing with the case $N\in\N$ is \cite{delespomu}, where, for any fixed positive integer  $N$,   it is proved the existence of  a solution to \eqref{eq} with $a=1$ and $\delta_0$ replaced by $\delta_{p_\la}$ for  a suitable $p_\la\in \Omega$   with $N+1$ blowing up points  at the vertices of a sufficiently tiny regular polygon  centered in $p_\la$; moreover $p_\la$ lies uniformly away from the boundary $\partial \Omega$  but its location  is determined by the geometry of the domain in an $\la-$dependent way and does not seem possible to be prescribed arbitrarily as in \cite{espo}.

 The case $N\in\N$ is more difficult to treat, and at the same time the most relevant to physical applications. Indeed, in vortex theory the number
$N$
represents vortex multiplicity, so that  in that context the most interesting case is
precisely that in which it is a positive integer. The difference between the case $N\in \N$ and $N\not\in \N$ is analitically essential. Indeed,  as usual in  
problems involving small parameters and concentration phenomena like \eqref{eq},  after
suitable rescaling of the blowing-up around a concentration point one sees a limiting equation.  
More specifically, as we will see in Section 2, we can associate to \eqref{eq} the limiting problem of Liouville type \eqref{limit} which will play a crucial role in the construction of  solutions blowing up at $0$ as $\la\to 0^+$;   
if $N\in\N$, \eqref{limit} admits a larger class of finite mass  solutions  with respect to the case $N\not\in\N$ since the family of all solutions extends to one carrying an extra parameter $b\in\R^2$  (see \cite{Pratar}). 

In this paper we are interested in  finding conditions on the potential $a$ and on the domain $\Omega$ under which there exists a solution $u_\la$ blowing up at $0$.  
Even though finding general conditions is a notoriously open
issue, our analysis reveals that  the interplay between the geometry of the domain, which is described in terms of the Robin function $H(x,x)$, and the potential $a$ plays a crucial role.
More specifically our conditions involve the first and the second derivative of $a$ and $H(x,x)$. 
\

Now we pass to provide the exact formulations of our results. 
In the following  we will assume that $$a\in C(\overline\Omega)\cap C^2(\Omega)\;\;\hbox{ and }\;\;\inf_\Omega a(x)>0.$$ Moreover,  after suitably
rotating the coordinate system, we may assume that in a small neighborhood of $0$ the following expansion holds:
$$a(x)=a(0)+\langle\nabla a(0), x\rangle +\frac{a_{11}x_1^2+a_{22}x_2^2}{2}+o(|x|^2)\hbox{ as}  x\to 0,$$
where $a_{ii}=\frac{\partial^2a}{\partial x_i^2}(0)$. 

\begin{thm}\label{th1} Let  $N\geq 2$, $N\in\N$. Assume that\footnote{Here $\nabla_xH(0,0)$ denotes the gradient of the function $x\mapsto H(x,0)$ at $0$.}  
\beq\label{assurd0}\nabla a(0)+4\pi(N+2)a(0)\nabla_xH(0,0)=0, \quad  \Delta a(0)\neq 16\pi^2(N+2)^2a(0)|\nabla_xH(0,0)|^2. \eeq
Then, for $\la$ sufficiently small, 
the problem \eqref{eq} has a family of solutions $u_\la$  blowing up at $0$ as $\la\to 0^+$. More precisely the following holds:  \beq\label{the1}\la a(x)e^{u_\la} dx \to 8\pi (1+N)\delta_0\eeq in the measure sense. More precisely $u_\la$ satisfies 
 \beq\label{the3}u_\la= 
 4\pi (N+2)G(x,0)+o(1)\eeq away from $0$. 
\end{thm}
In the particular case when $0$ is a critical point of the potential $a$ and of the Robin function,  we get the existence of a solution blowing up at $0$ provided that $\Delta a(0)\neq 0$. 

\begin{cor}\label{th1cor} Let  $N\geq 2$, $N\in\N$.  Assume that
$$\nabla a(0)=\nabla_xH(0,0)=0, \quad  \Delta a(0)\neq 0. $$
Then, for $\la$ sufficiently small 
the problem \eqref{eq} has a family of solutions $u_\la$  blowing up at $0$ as $\la\to 0^+$.  More precisely $u_\la$ satisfies \eqref{the1}-\eqref{the3} of Theorem \ref{th1}. 
\end{cor}

The case $N=1$ is considered in  a separate theorem since the result requires different assumptions; the above result continues to hold in symmetric domains under an additional  relation involving the second derivatives of $a$.

\begin{thm}\label{th2} Let $N=1$. Assume that  $\Omega$ is $\ell$-symmetric for some $\ell\geq3$, i.e. \beq\label{assurd} x\in\Omega\Longleftrightarrow e^{{\rm i}\frac{2\pi}{\ell}} x\in\Omega\eeq 
and, in addition, \beq\label{assurd2}\nabla a(0)
=0, \quad  \Delta a(0)\neq 0, \quad a_{11}=a_{22}.\eeq 
Then, for $\la$ sufficiently small 
the problem \eqref{eq} has a family of solutions $u_\la$  blowing up at $0$ as $\la\to 0^+$. More precisely $u_\la$ satisfies \eqref{the1}-\eqref{the3} of Theorem \ref{th1}. 
\end{thm}

We point out that  in symmetric domains  the center of symmetry $0$ is a critical point of the Robin function, so the condition $\nabla_xH(0,0)=0$ is automatically satisfied. Assumptions \eqref{assurd}-\eqref{assurd2} are obviously satisfied if $\Omega$ is a ball centered at $0$ and $a$ is a radially symmetric potential with a nondegenerate critical point at $0$.

\

The proofs use singular perturbation methods. Roughly speaking, the first step consists in the construction of an approximate solution, which should turn out to  be precise enough. In view of the expected asymptotic behavior,  the shape of such  approximate solution will resemble a \textit{bubble} of the form \eqref{bubble} with  a suitable choice of the parameter $\de=\de(\la,b)$. Then  we look for a solution to \eqref{eq} in a small neighborhood of the first approximation.  As quite standard in singular perturbation theory, a crucial ingredient is nondegeneracy of
the explicit family of 
solutions of the limiting Liouville problem \eqref{limit}, 
in the sense that all bounded elements in the kernel of the linearization
correspond to variations along the parameters of the family, as established in \cite{delespomu}.
This allows us to study the invertibility of the linearized operator associated to the problem \eqref{eq}  under suitable orthogonality conditions. Next we introduce an intermediate problem and  a fixed point argument will provide a solution for an  auxiliary equation, which turns out to be solvable for any choice of $b$. Finally we test the auxiliary equation on the elements of the kernel of the linearized operator  and we find out that, 
 in order to find an \textit{exact} solution of \eqref{eq}, the parameter  $b$ should be a zero for 
 a \textit{reduced} finite dimensional map. 

\

The rest of the paper is organized as follows. Section 2 is devoted to some preliminary
results, notation, and the definition of the approximating solution. Moreover, a more
general version of Theorems \ref{th1}-\ref{th2} is stated there (see Theorems \ref{main1}-\ref{main2}). The error up to which the approximating solution solves problem \eqref{eq} is estimated in Section 3.   In Section 4  we
prove the solvability of the linearized problem.  Section 5 considers the solvability of an auxiliary problem by a  contraction argument.. Finally, in Section 6, we prove the existence
results and we conclude the proof of Theorems \ref{th1}-\ref{th2}. In Appendix A and Appendix B we collect some results, most of them well-known, 
 which are usually referred to throughout the
paper.

\section{Preliminaries and statement of the main results}
We are going to provide an equivalent formulation of problem \eqref{eq} and Theorems \ref{th1}-\ref{th2}. Indeed,  let us observe that, setting  $v$ the regular part of $u$, namely \beq\label{chva}v= u+4\pi (\alpha-1) G(x, 0),\quad \alpha=N+1,\eeq  problem \eqref{eq} is then equivalent to solving the following (regular) boundary value problem
\beq\label{proreg}\left\{\begin{aligned} & -\Delta v=\la V(x)|x|^{2(\alpha-1)}e^v&\hbox{ in }&\Omega\\ &v=0&\hbox{ on }&\partial \Omega\end{aligned}\right.,\eeq
where $V(x)$ is the new potential \beq\label{aaa}V(x)=a(x)e^{-4\pi  (\alpha-1) H(x,0)} .\eeq
Here $G$ and $H$ are Green's function and its regular part as defined in the introduction. 
This problem is actually variational. Indeed, let us consider the following energy functional associated with \eqref{proreg}:
$$J(v)=\frac12\into |\nabla v|^2 dx-\la \into V(x)|x|^{2(\alpha-1)}e^vdx,\quad v\in H^1_0(\Omega).$$
Then the following Moser-Trudinger inequality 
(\cite{Moe, Tru}) guarantees that $J$ is of class $C^1(H^1_0(\Omega))$ and solutions of \eqref{proreg} correspond to critical points of $J$.

\begin{lemma}\label{tmt} There exists $C>0$ such that for any bounded domain $\Omega$ in $\rr^2$
 $$\into e^{\frac{4\pi u^2}{\|u\|^2}}dx\le C |\Omega|\quad \forall u\in{ H}^1_0(\Omega),$$ where  $|\Omega|$ stands for the measure of the domain $\Omega$.
 In particular,  for any $q\geq 1$
 $$\| e^{u}\|_{q}\le  C^{\frac1q} |\Omega|^{\frac1q} e^{{q\over 16\pi}\|u\|^2}\quad \forall  u\in{H}^1_0(\Omega).$$

\end{lemma}

Theorems \ref{th1}-\ref{th2} will be a consequence of  more general results concerning Liouville-type problem \eqref{proreg}. In order to provide such  results  \eqref{proreg}, we now give a construction of a suitable approximate solution for \eqref{proreg}. In what follows, we identify $x=(  x_1,x_2)\in \R^2$ with $x_1+{\rm i}x_2\in \C$.  Moreover, $\langle x_1, x_2\rangle$ stands for the inner product between the vectors $x_1, x_2\in\R^2$, whereas $x_1 x_2$ will denote the multiplication of the complex numbers $x_1$, $x_2$. Clearly $\langle x_1,  x_2\rangle=\Re(x_1 \overline x_2)$.

For any $\al\in\N$, we can associate to \eqref{proreg} a limiting problem of Liouville type which will play a crucial role in the construction of the blowing-up solutions: 
\beq\label{limit}
-\Delta w=|x|^{2(\alpha-1)}e^w\quad \hbox{in}\;\; \rr^2,\qquad
\int_{\R^2} |x|^{2(\al-1)}e^{w(x)}dx<+\infty.
\eeq
A complete classification for solutions of \eqref{limit} is due to \cite{Pratar} and corresponds, 
 in complex notation,  to the three-parameter family of functions
\beq\label{bubble}
w^\al_{\de,b}(x):=\log  {8\al^2\de^{2\al}\over (\de^{2\al}+|x^{\al}- b|^2)^2}\quad
\de>0,\,b\in \C.
\eeq
The following quantization property holds: \beq \label{quantum} \int_{\R^2}
|x|^{2(\al-1)}e^{w_{\de,b}^\al(x)}dx = 8 \pi \alpha .\eeq
 In the following we  agree that 
$$W_\la=w^{\al}_{\de,b}(x),$$ where 
the value $\delta=\delta(\la,b)$ is defined as:

\begin{equation} \label{delta} \delta^{{2\alpha}} :=\frac{\la}{8\al^2}V(0)e^{8\pi\sum_{i=1}^\al H(0,\beta_i)}.  \end{equation}

To obtain a better first approximation, we need to modify the functions $W_\la$   in order to satisfy the zero boundary condition. Precisely, we consider the projections $P W_\la $ onto the space $ H^1_0(\Omega)$ of
$W_\la$, where the projection  $P:H^1(\R^N)\to  H^1_0(\Omega)$ is
defined as the unique solution of the problem
$$
 \Delta P v=\Delta v\quad \hbox{in}\ \Omega,\qquad  P v=0\quad \hbox{on}\ \partial\Omega.
$$
Let us consider $b$ in a small neighborhood of $0$ and let us denote by $\beta_0,\ldots,  \beta_{\alpha-1}$ the $\alpha$-roots of $b$, i.e., $\beta_i^\alpha=b$ and $\beta_i\neq\beta_h$ for $i\neq h$. Observe that the function  $\sum_{i=0}^{\alpha-1}  H(x, \beta_i)$ is harmonic in $\Omega$ and satisfies $\sum_{i=0}^{\alpha-1}  H(x, \beta_i)=\frac{1}{2\pi}\log|x^\alpha-b|$ on $\partial \Omega.$ A straightforward computation gives that for any $x\in\partial\Omega$ 
$$\bigg|PW_\la- W_\la+\log\(8\al^2\de^{2\al}\)-8\pi\sum_{i=0}^{\alpha-1}  H(x, \beta_i)\bigg|=\bigg|W_\la-\log\(8\al^2\de^{2\al}\)+4 \log|x^\alpha-b|
\bigg|\leq C\de^{2\alpha}.$$
Since the expressions considered inside the absolute values 
are harmonic in $\Omega$, then the maximum principle applies and implies
the following asymptotic expansion
\beq\label{pro-exp1}\begin{aligned}
 PW_\la=& W_\la-\log\(8\al^2\de^{2\al}\)+8\pi\sum_{i=0}^{\alpha-1} H(x, \beta_i)+O(\de^{2\al})\\ =&-2\log\(\de^{{2\alpha}}+|x^{\al}- b|^2\)+8\pi\sum_{i=0}^{\alpha-1} H(x, \beta_i)+O(\de^{2\al})
\end{aligned}\eeq 

uniformly for $x\in \bar\Omega$ and $b $ in a small neighborhood of $0$.

We shall look for a solution to \eqref{proreg} in a small neighborhood of the first approximation, namely a solution of the form
 $$v_\la=PW_\la
 + {\phi}_\la,$$ where the rest term
$\phi_\la$ is small in
$H^1_0(\Omega)$-norm.

We are now in the position to state  the main theorems of the paper.

\begin{thm} \label{main1} Assume that $\alpha\geq 3$ and hypotheses \eqref{assurd0}  hold. Then, for $\la$ sufficiently small 
 there
exist  $\phi_\la \in H^1_0(\Omega)$ and $b=b_\la=O(\la^{\frac{\al+1}{2\al}})$ such that the couple $
PW_\la+\phi_\la$ solves problem
\eqref{proreg}. 
Moreover, for any fixed $\e>0$,
\beq\label{barb}\| \phi_\la \|_{H^1_0(\Omega)}\leq \la^{\frac{1}{\alpha}-\e}\;\;\hbox{ for } \la \hbox{ small enough}.\eeq
\end{thm}

\begin{thm} \label{main2} Assume that $\al=2$, and hypotheses \eqref{assurd}-\eqref{assurd2}  hold. Then, for $\la$ sufficiently small there
exist 
$\phi_\la \in H^1_0(\Omega)$ and $b=b_\la=O(\la^{\frac{\al+1}{2\al}})$ such that the couple $
PW_\la+\phi_\la$ solves problem
\eqref{proreg}. 
Moreover, for any fixed $\e>0$, \eqref{barb} holds. 
\end{thm}

In the remaining part of the paper we will prove Theorems \ref{main1}-\ref{main2} and at the end of the Section 6 we shall see how Theorems \ref{th1}-\ref{th2} follow quite directly as a corollary according to \eqref{chva} and  \eqref{aaa}.

We end up this section by setting the notation and basic well-known
facts which will be of use in the rest of the paper. We  denote by  $\|\cdot\|$ and $\|\cdot\|_p$  the norms in  the space $ H^1_0(\Omega)$ and $L^p(\Omega)$, respectively, namely
 \beq\label{nott}\|u\|:=\|u\|_{ H^1_0(\Omega)}
 ,\qquad \|u\|_p:=\|u\|_{L^p(\Omega)}
 \quad \forall u\in  H^1_0(\Omega).\eeq

For any
$\alpha\ge1$ we will make use of the Hilbert spaces
\begin{equation}\label{ljs}
\mathrm{L}_\alpha (\rr^2):=\left\{u \in L^2_{loc}(\rr^2)\ :\  \left\|{|y|^{\alpha-1} \over
1+|y|^{{2\alpha}}}u\right\|_{{L}^2(\rr^2)}<+\infty\right\}\end{equation}
 and
\begin{equation}\label{hjs}\mathrm{H}_\alpha (\rr^2):=\left\{u\in {\rm W}^{1,2}_{loc}(\rr^2) \ :\ \|\nabla u\|_{{L}^2(\rr^2)}+\left\|{|y|^{\al-1} \over 1+|y|^{{2\alpha}}}u\right\|_{{L}^2(\rr^2)}<+\infty\right\},\end{equation}
 endowed with the norms
$$\|u\|_{\mathrm{L}_\alpha }:= \left\|{|y|^{\al-1} \over 1+|y|^{2\alpha}}u\right\|_{{L}^2(\rr^2)}\
\hbox{and }\ \|u\|_{\mathrm{H}_\alpha }:= \(\|\nabla
u\|^2_{{L}^2(\rr^2)}+\left\|{|y|^{\al-1} \over
1+|y|^{{2\alpha}}}u\right\|^2_{{L}^2(\rr^2)}\)^{1/2}.$$
We  denote by $\langle u,v \rangle_{\mathrm{L}_\alpha }$ and $\langle u,v \rangle_{\mathrm{H}_\alpha }$
the natural scalar product in ${\mathrm{L}_\alpha }$ and in ${\mathrm{H}_\alpha }$, respectively.

\begin{prop}\label{compact}
The embedding $i_\al:\mathrm{H}_\alpha
(\rr^2)\hookrightarrow\mathrm{L}_\alpha (\rr^2)$ is compact.
\end{prop}
\begin{proof}
  See \cite[Proposition 6.1]{gpistoia}.
\end{proof}

As commented in the introduction, our proof uses the singular
perturbation methods. For that, the nondegeneracy of the functions
that we use to build our approximating solution is essential. Next
proposition is devoted to the nondegeneracy of the finite mass 
solutions of the Liouville equation (regular and singular).

\begin{prop}
\label{esposito} Assume that $\phi:\R^2\to\R$  solves the problem
\begin{equation}\label{l1}
-\Delta \phi =8\alpha^2{|y|^{2(\alpha-1)}\over (1+|y^{\al}-\xi|^2)^2}\phi\;\;
\hbox{in}\ \rr^2,\quad \int_{\R^2}|\nabla
\phi(y)|^2dy<+\infty.
\end{equation}
 Then there exist $c_0,\,c_1,\, c_2\in\R$
such that
$$\phi(y)=c_0  Z_0+ c_1Z_1 +c_2Z_2.$$
$$Z_0(y):   = {1-|y^{\al}-\xi|^2\over 1+|y^{\al}-\xi|^2} ,\ \; \;Z_1(y):={ \Re(y^{\al}-\xi)\over  1+|y^{\al}-\xi|^2} 
,\ \;\;Z_2(y):={ \Im(y^{\al}-\xi)\over  1+|y^{\al}-\xi|^2}.
$$ \end{prop}
\begin{proof}
In \cite[Theorem 6.1]{gpistoia} it was proved that any solution
$\phi$ of \eqref{l1} is actually a bounded solution. Therefore we
can apply the result in \cite{dem} to conclude that $\phi= c_0 \phi_0 + c_1
\phi_1 + c_2 \phi_2$ for some $c_0,c_1,c_2\in \R$.

\end{proof}

In our estimates throughout the paper, we will frequently denote by $C>0$, $c>0$ fixed
constants, that may change from line to line, but are always
independent of the variable under consideration. We also use the
notations $O(1)$, $o(1)$, $O(\lambda)$, $o(\lambda)$ to describe
the asymptotic behaviors of quantities in a standard way.

\section{Estimate of the error term}
The goal of this section is to  provide an estimate of the error up to which the function $W_\la$ solves problem \eqref{proreg}.
\begin{lemma}\label{aux} Let $r>0$ be a fixed    number. Define
$${R}_\la:=
-\Delta PW_\la-\la V(x) |x|^{2(\alpha-1)}e^{PW_\la}= |x|^{2(\alpha-1)}e^{W_\la}-\la V(x)|x|^{2(\alpha-1)} e^{PW_\la}.$$
 For any fixed $p\geq 1$ the following holds
$$\|R_\la\|_{p}=O(\la^{\frac{1}{\alpha p}}).$$ uniformly for $|b|\leq r\sqrt{\la}$.  Consequently, for every fixed  $p\geq 1$,
\beq\label{judo}\|\la V(x)|x|^{2(\alpha-1)} e^{PW_\la}\|_p=\||x|^{2(\alpha-1)}e^{W_\la}\|_p+o(1)=O(\la^{\frac{1-p}{\alpha p}})\eeq
uniformly for $|b|\leq r\sqrt{\la}$.
\end{lemma}

\begin{proof}
By \eqref{pro-exp1}  and the choice of $\delta$ in \eqref{delta} we derive  \beq\label{cuore}\begin{aligned}&\la V(x)|x|^{2(\alpha-1)} e^{PW_\la}\\ &
=\frac{\la}{8\al^2\de^{2\al}} V(x)|x|^{2(\alpha-1)} e^{W_\la+8\pi \sum_{i=0}^{\al-1}H(x, \beta_i) +O(\de^{2\alpha} )} 
\\ &=|x|^{2(\alpha-1)}e^{W_\la}\frac{V(x)}{V(0)}e^{8\pi\sum_{i=0}^{\alpha-1}(H(x, \beta_i)-H(0,\beta_i))+O(\de^{2\alpha} )}
\\ &=|x|^{2(\alpha-1)}e^{W_\la}\frac{a(x)}{a(0)}e^{-4\pi(\al-1) (H(x,0)-H(0,0))+8\pi\sum_{i=0}^{\al-1}(H(x, \beta_i)-H(0,\beta_i))+O(\de^{2\alpha} )}
.\end{aligned}\eeq According to \eqref{expa} we have 
$$H(x,0)-H(0,0)=
\Re\Big(\frac{d\tilde H}{dx} (0,0) x\Big)+O(|x|^{2})
$$
whereas, by   Lemma \ref{robin10},
$$\sum_{i=0}^{\al-1}(H(x, \beta_i)-H(0,\beta_i))=
\al \Re\Big(\frac{d\tilde H}{dx} (0,0) x\Big)+O(|x|^{2})+O(|b|^2)
$$ by which we arrive at 
$$\begin{aligned}&\la V(x)|x|^{2(\alpha-1)} e^{PW_\la}\\&
=|x|^{2(\alpha-1)}e^{W_\la}\frac{a(x)}{a(0)}
e^{4\pi(\al+1)\Re(\frac{d\tilde H}{dx} (0,0) x)+O(|x|^{2})+O(|b|^2)+O(\de^{2\alpha} )} \\ &=|x|^{2(\alpha-1)}e^{W_\la}\frac{a(x)}{a(0)}
\Big( 1+ 4\pi(\al+1)\Re\Big(\frac{d\tilde H}{dx} (0,0) x\Big)+ O(|x|^{2})+O(|b|^2)+O(\de^{2\alpha} )\Big)
\\ &=|x|^{2(\alpha-1)}
e^{W_\la}\Bigg(1+
\frac{\langle\nabla a(0), x\rangle}{a(0)} + 4\pi(\al+1)\Re\Big(\frac{d\tilde H}{dx} (0,0) x\Big)+ O(|x|^{2})+O(|b|^2)+O(\de^{2\alpha} )\Bigg) 
.
\end{aligned}
$$ 
Taking into account that $\Re(\frac{d\tilde H}{dx} (0,0) x)=\langle \nabla_x H(0,0),x\rangle$ and using that 
$$\frac{\langle \nabla a(0), x\rangle}{a(0)} + 4\pi(\al+1)\Re\Big(\frac{d\tilde H}{dx} (0,0) x\Big)=\frac{\langle\nabla a(0) ,x\rangle}{a(0)}+ 4\pi(\al+1)\langle  \nabla_x H(0,0),x\rangle=0$$ by assumptions \eqref{assurd0} and \eqref{assurd}-\eqref{assurd2},
we arrive at
$$\begin{aligned}&\la V(x)|x|^{2(\alpha-1)} e^{PW_\la}=|x|^{2(\alpha-1)}
e^{W_\la}+ \big(O(|x|^{2})+O(|b|^2)
+O(\de^{2\alpha} )\big)|x|^{2(\alpha-1)}e^{W_\la}.\end{aligned}$$
Now if we scale $x=\de y$, recalling that $|b|\leq \sqrt{\la}\leq C\de^\al$,  we get 
$$\begin{aligned}
|x|^{2(\alpha-1)}
e^{W_\la}&
=8\al^2
\frac{|y|^{2(\al-1)}}{\de^2(1+|y^{\alpha }-\de^{-\al}b|^2)^2}
= O\bigg(\frac{1}{\de^2(1+|y|^{2\al+2})}\bigg)
.\end{aligned} $$
and, similarly,
$$\begin{aligned}|x|^{2\alpha}
e^{W_\la}&
=8\al^2
\frac{|y|^{2\al}}{(1+|y^{\alpha }-\de^{-\al}b|^2)^2}
= O\bigg(\frac{1}{1+|y|^{2\al}}\bigg)
.\end{aligned} $$
by which
$$\||x|^{2(\alpha-1)}
e^{W_\la}\|_p=O(\de^{2\frac{1-p}{p}}) \qquad 
\||x|^{2\alpha}
e^{W_\la}\|_p=O(\de^{\frac2p})$$
The thesis is thus proved.

\end{proof}

\section{Analysis of the linearized operator}
According to Proposition \ref{esposito}, by the change of variable $x=\de y$, we immediately get that  all solutions $\psi \in \mathrm{H}_\al(\R^2)$  of
$$
-\Delta \psi= 8{\alpha}^2{\de^{2\alpha}|x|^{2(\alpha-1)}\over (\de^{2\alpha}+|x^{\alpha }-b|^2)^2}\psi =|x|^{2(\alpha-1)} e^{W_\la}\psi\quad \hbox{in}\quad \rr^2$$
are linear combinations of the functions
$$Z^0_{\delta,b}(x)={\delta^{2\alpha}-|x^{\alpha }-b|^2\over \delta^{2\alpha}+|x^{\alpha }-b|^2},\ Z^1_{\delta,b}(x)=
{ \de^{\alpha }\Re(x^{\alpha }-b)\over  \de^{2\alpha}+|x^{\alpha }-b|^2},\ Z^2_{\delta,b}(x)=
{ \de^{\alpha }\Im(x^{\alpha }-b)\over  \de^{2\alpha}+|x^{\alpha }-b|^2}
.$$
We introduce their projections $PZ^j_{\delta,b}$ onto $H^1_0(\Omega).$ It is immediate that
 \begin{equation}\label{pz0}
PZ^0_{\delta,b}(x)=Z^0_{\delta,b}(x)+1+ O\(\de^{2\alpha}\)
\end{equation}
 and
  \begin{equation}\label{pzi}
PZ^j_{\delta,b}(x)=Z^j _{\delta,b}(x) + O\(\de^{\alpha }\),\;\; j=1,2
\end{equation}
uniformly with respect to $x\in\overline\Omega$ and $b$ in a small neighborhood of $0$.
\\
We agree that $Z_\la^j:=Z_{\delta,b}^j$ for any $j=0,1,2$, where $\delta$ is defined in terms of  $\la$ and $b$ according to \eqref{delta}.

Let us consider the following linear problem: given
$ h\in H_0^1(\Omega)$,  find a function $\phi\in  H^1_0(\Omega)$ satisfying
\begin{equation}\label{lla}
\left\{\begin{aligned}&-\Delta \phi   -\la V(x)|x|^{2(\al-1)} e^{P{W}_\la}\phi=\Delta h\\ &\into \nabla \phi\nabla PZ_\la^j=0\;\;j=1,2\end{aligned}\right..
\end{equation}

Before going on, we recall the following identities which follow by straightforward computations using Lemma \ref{copy}: for every $\xi\in \R^2$

\beq\label{id1}\begin{aligned}\intr |y|^{2(\alpha-1)}\log (1+|y^{\alpha }-\xi|^2)\frac{1-|y^{\alpha }-\xi|^{2}}{(1+|y^{\alpha }-\xi|^2)^3}dy&=\frac{1}{{\alpha}}\intr \log (1+|y|^2)\frac{1-|y|^{2}}{(1+|y|^{2})^3}dy\\ &=-\frac{\pi}{{2\alpha}},\end{aligned}\eeq
\beq\label{id2}\intr\frac{|y|^{2(\alpha-1)}}{(1+|y^{\alpha }-\xi|^2)^2}\frac{1-|y^{\alpha }-\xi|^{2}}{1+|y^{\alpha }-\xi|^{2}}dy=\frac{1}{{\alpha}}\intr\frac{1-|y|^{2}}{(1+|y|^{2})^3}dy=0,\eeq

\beq\label{id3}\begin{aligned}\intr\frac{|y|^{2(\alpha-1)}(\Re (y^{\alpha }-\xi))^2}{(1+|y^{\alpha }-\xi|^{2})^4}dy&=\intr\frac{|y|^{2(\alpha-1)}(\Im (y^{\alpha }-b))^2}{(1+|y^{\alpha }-\xi|^{2})^4}dy\\ &=\frac{1}{{2\alpha}}\intr\frac{|y|^2}{(1+|y|^{2})^4}dy=\frac{\pi}{12{\alpha}}.\end{aligned}\eeq

\begin{prop}\label{inv} Let $r>0$ be fixed. 
There exist $\lambda_0>0$ and $C>0$ such that for any $\la \in(0, \la_0)$, any $b\in \R^2$ with $|b|<r\sqrt{\la}$ and   any $h\in  H^1_0(\Omega)$, if  
$\phi\in  H^1_0(\Omega)$ solves  \eqref{lla}, then the following holds $$\|\phi\|    \leq C  |\log\la |    \|h\| .$$
\end{prop}
\begin{proof}
We argue by contradiction. Assume that there exist sequences
$\la_n\to0,$ $ h_n\in  H^1_0(\Omega)$, $|b_n|\leq r\sqrt{\la_n}$ and $\phi_n\in  H^1_0(\Omega)$ 
which solve \eqref{lla} and
\begin{equation}\label{inv2}
\|\phi_n\| =1, \qquad
   |\log\la_n | \| h_n\| \to 0.\end{equation}
Let $\de_n>0$ be the value associated to $\la_n$ according to  \eqref{delta}. Then we may assume $$\de_n^{-\al}b_n\to b_0.$$ We define  $\widetilde {\Omega}_n :={\Omega\over {\de}_n }$ and
 $$  \tilde \phi_n(y):=\left\{\begin{aligned}&\phi_n\({\de}_n y\)&\hbox{ if }&y\in \widetilde{\Omega}_n \\ &0&\hbox{ if }&y\in \rr^2\setminus\widetilde{\Omega}_n \end{aligned} \right. . $$
\\
 In what follows at many steps of the arguments we will pass to a subsequence, without further notice.  
We split the remaining argument into five steps.
\medskip

\noindent{\em Step 1. We will show that
\begin{equation*}\label{step1.0} \tilde\phi_n \;\; \hbox{ is bounded in }\mathrm{H}_{\al} (\rr^2).
\end{equation*}
  }

It is immediate to check that
\beq\label{bounabla}\int_{\R^2}|\nabla \tilde\phi_n|^2dy=\into|\nabla  \phi_n|^2dx\le1.\eeq
Next, we multiply the  equation  in \eqref{lla} by $\phi_n$; then  we integrate over $\Omega$  to obtain

$$\begin{aligned}\la_n \into  V(x)|x|^{2(\al-1)}e^{PW_{\la_n}}\phi_n^2dx=& \into |\nabla\phi_n|^2dx +\into \nabla  h_n\nabla \phi_ndx\end{aligned}$$
which implies, by \eqref{inv2},
\beq\label{W22}\la_n \int_{\Omega} V(x) |x|^{2(\al-1)}e^{PW_{\la_n}}\phi_n^2dx \leq
C.\eeq So, Lemma \ref{aux}  gives $\into |x|^{2(\al-1)}
e^{W_{\la_n}}\phi_n^2\leq C$ or,
equivalently,
$$\int_{\R^2} {|y|^{2(\al-1)}\over \(1+|y^{\alpha }-\de_n^{-\al}b_n|^2\)^2} \tilde\phi_n^2  dy \le C.$$
Combining this with \eqref{bounabla}, we deduce that    $\tilde \phi_n$ is bounded in the space $\mathrm{H}_{\al} (\rr^2)$.

\bigskip

\noindent{\em Step 2. We will show that, for some $ \gamma_0 \in\rr $,
$$\tilde\phi_n\to \gamma_0  \frac{1-|y^{\alpha }-b_0|^{2}}{1+|y^{\alpha }-b_0|^{2}}\;\ \hbox{  weakly in }\rm{H}_{\al} (\rr^2) \hbox{ and strongly in }\rm{L}_{\al} (\rr^2) .$$
  }

Step 1 and Proposition \ref{compact} give
 $$\tilde \phi_n\to f\;\hbox{ weakly in }\mathrm{H}_\al (\rr^2) \hbox{ and strongly in }\mathrm{L}_\al(\rr^2) .$$ Let $\tilde\psi\in C^\infty_c(\R^2)$ and set
 ${\psi}_n=\tilde{\psi}(\frac{x}{{\de}_n})\in C^\infty_c(\Omega)$, for large $n$.
 We multiply the  equation in \eqref{lla} by ${\psi}_n,$   we integrate over $\Omega$ and we get
\begin{equation}\label{1.1.1.1}\int_{\widetilde{\Omega}_n} \nabla{\tilde\phi_n}\nabla\tilde\psi dy - \la_n\int_{{\Omega}} V(x)|x|^{2(\al-1)} e^{PW_{\la_n}}\phi_n \psi_n dx=-\int_{{\Omega}}\nabla
 h_n\nabla {\psi}_ndx.\end{equation}
 According to Lemma \ref{aux} we have 
 $$\begin{aligned}\la_n\int_{{\Omega}}V(x) |x|^{2(\al-1)} e^{PW_{\la_n}}\phi_n \psi_n dx&=\int_{{\Omega}} |x|^{2(\al-1)} e^{W_{\la_n}}\phi_n \psi_n dx+o(1)
 \\ &=8{\alpha}^2\int_{\R^2} {|y|^{2(\al-1)}\over (1+|y^{\alpha }-\de_n^{-\al}b_n|^2)^2} \tilde\phi_n \tilde\psi dy+o(1)\\ &=8{\alpha}^2\int_{\R^2} {|y|^{2(\al-1)}\over (1+|y^{\alpha }-b_0|^2)^2} f \tilde\psi dy+o(1).\end{aligned}$$
 Finally, by \eqref{inv2}, using that $\into |\nabla {\psi}_n|^2=\int_{\R^2}|\nabla\tilde \psi|^2$,
\beq\label{macro}\int_{{\Omega}}|\nabla
{h}_n\nabla {\psi}_n|dx=O(\|{h}_n\|)=o(1).\eeq
\medskip Therefore, we may pass to the limit in \eqref{1.1.1.1} to
obtain
\begin{align}\nonumber&\int\limits_{\R^2} \nabla f \nabla\tilde\psi dy  =8{\alpha}^2\int_{\R^2} {|y|^{2(\al-1)}\over (1+|y^{\alpha }-b_0|^2)^2} f \tilde\psi dy
.\end{align}
Thus 
 $f$ 
solves the equation
$$-\Delta f=8{\alpha}^2 {|y|^{2(\al-1)}\over (1+|y^{\alpha }-b_0|^2)^2} f .$$
Proposition \ref{esposito} gives  \beq\label{crucial3}f=\gamma_0Z_0+\gamma_1Z_1
+\gamma_2Z_2\eeq
for some $\gamma_0,\,\gamma_1,\, \gamma_2\in \R$. It remains to show that $\gamma_1=\gamma_2=0.$ Indeed, we compute $$\begin{aligned}0&=\into \nabla \phi_n\nabla PZ_{\la_n}^1dx=\into|x|^{2(\al-1)} e^{W_{\la_n}} \phi_n Z_{\la_n}^1dx=8{\alpha}^2\intr {|y|^{2(\alpha-1)}\over (1+|y^{\alpha }-\de_n^{-\alpha}b_n|^2)^2}\tilde\phi_n Z_1dy\\ &=
8{\alpha}^2\int_{\R^2}f {|y|^{2(\alpha-1)}\over (1+|y^{\alpha }-b_0|^2)^2}Z_1 dy+o(1)=
\int_{\R^2}\nabla f\nabla Z_1dy+o(1). \end{aligned}$$ We get $\int_{\R^2} \nabla f \nabla Z_1=0$, by which, taking into account  that $\int_{\R^2} \nabla Z_1\nabla Z_0=\int_{\R^2} \nabla Z_1\nabla Z_2=0$, 
$$\gamma_1\int_{\R^2}| \nabla PZ_1|^2dy=0. $$
So $\gamma_1=0$ and, similarly, $\gamma_2=0.$

\bigskip

\noindent{\em Step 3. We will show that $$\intr\frac{|y|^{2(\alpha-1)}}{(1+|y^{\alpha }-\de_n^{-\al}b_n|^2)^2}\tilde\phi_n dy=o\Big(\frac{1}{\log\la_n}\Big).$$}

We multiply the equation in \eqref{lla} by $PZ^0_{\la_n}$, we integrate over $\Omega$ and we get 
\beq\label{pi0}\into \nabla \phi_n\nabla PZ^0_{\la_n}dx-\la_n\into V(x)|x|^{2(\alpha-1)}e^{PW_{\la_n}}\phi_nPZ^0_{\la_n} dx=-\into \nabla h_n \nabla PZ^0_{\la_n} dx.\eeq
We are now concerned with the estimates of each term of the above expression. 
First, we compute
\beq\label{pi1}\into \nabla \phi_n\nabla PZ^0_{\la_n}dx=
\into |x|^{2(\alpha-1)}e^{W_{\la_n}} \phi_n Z^0_{\la_n}dx.\eeq
Using Lemma \ref{aux} (with $p=2$) and \eqref{pz0}, we obtain
\beq\label{pi2}\begin{aligned}\la_n\into& V(x)|x|^{2(\alpha-1)}e^{PW_{\la_n}}\phi_nPZ^0_{\la_n} dx=\into |x|^{2(\alpha-1)}e^{W_{\la_n}}\phi_n(Z^0_{\la_n}+1) dx+O(\la_n^{\frac{1}{{2\alpha }}})
\\ &=\into |x|^{2(\alpha-1)}e^{W_{\la_n}}\phi_nZ^0_{\la_n} dx+8{\alpha}^2\intr\frac{|y|^{2(\alpha-1)}}{(1+|y^{\alpha }-\de_n^{-\al}b_n|^2)^2}\tilde\phi_n dy+O(\la_n^{\frac{1}{{2\alpha }}}) .\end{aligned}\eeq
Finally, since $PZ^0_{\la}=O(1)$, we have $\into |\nabla PZ^0_{\la}|^2=\into |x|^{2(\alpha-1)}e^{W_{\la}}PZ^0_{\la}=
O(1)$, by which, owing to \eqref{inv2},
\beq\label{pi3} \into |\nabla h_n|\,| \nabla PZ^0_{\la_n}| dx\leq \|h_n\|\, \|PZ^0_{\la_n}\|=o\Big(\frac{1}{\log\la_n}\Big).\eeq
We now multiply \eqref{pi0} by $\log\la_n$ and pass to the limit: inserting \eqref{pi1}, \eqref{pi2}, \eqref{pi3}, we obtain the thesis of the step.

\bigskip
\noindent{\em Step 4. We will show that $\gamma_0=0$.}

\medskip

 We multiply the equation in \eqref{lla} by $PW_{\la_n}$, we integrate over $\Omega$ and we get 
\beq\label{ti0}\into \nabla \phi_n\nabla PW_{\la_n}dx-\la_n\into V(x) |x|^{2(\alpha-1)}e^{PW_{\la_n}}\phi_nPW_{\la_n} dx=-\into \nabla h_n \nabla PW_{\la_n}dx.\eeq
Let us estimate each of the terms above. Let us begin with:
\beq\label{ti1}\into \nabla \phi_n\nabla PW_{\la_n}dx=\into |x|^{2(\alpha-1)} e^{W_{\la_n}}\phi_n dx=8\alpha^2\intr \frac{|y|^{2(\alpha-1)}}{(1+|y^{\alpha }-\de_n^{-\al}b_n|^2)^2}\tilde\phi_ndy=o(1)
\eeq by step 3. 
By Lemma \ref{aux} and  \eqref{inv2}, using that $|PW_{\la_n}|=O(|\log\la_n|)$, we get 
\beq\label{ti2}\begin{aligned}\la_n\into V(x)|x|^{2(\alpha-1)}e^{PW_{\la_n}}&\phi_nPW_{\la_n} dx=\into |x|^{2(\alpha-1)}e^{W_{\la_n}}\phi_nPW_{\la_n} dx+o(1)
\\ &=8{\alpha}^2\intr \frac{|y|^{2(\alpha-1)}}{(1+|y^{\alpha }-\de_n^{-\al}b_n|^2)^2}\tilde\phi_nPW_{\la_n} (\de_n y)dy+o(1). \end{aligned}
\eeq
Observe that by \eqref{pro-exp1} $$PW_{\la_n}(\de_n y)=-2\log (1+|y^{\alpha }-\de_n^{-\al}b_n|^2)+8\pi {\alpha} H(\de_n y,0)-4{\alpha}\log\de_n+O(\sqrt{\la_n})$$
by which 
$$PW_{\la_n}(\de_n y)+4{\alpha}\log\de_n\to -2\log (1+|y^{\alpha }-b_0|^2)+8\pi {\alpha} H(0,0)\;\;\hbox{ uniformly in }\R^2.$$
Using this convergence in \eqref{ti2}, and recalling  step 2, we obtain
$$\begin{aligned}\la_n\into &V(x)|x|^{2(\alpha-1)}e^{PW_{\la_n}}\phi_nPW_{\la_n} dx\\ &=-16{\alpha}^2\gamma_0\intr \log (1+|y^{\alpha }-b_0|^2)\frac{|y|^{2(\alpha-1)}}{(1+|y^{\alpha }-b_0|^2)^2}\frac{1-|y^{\alpha }-b_0|^{2}}{1+|y^{\alpha }-b_0|^{2}}dy \\ &\;\;\;\;+64\pi {\alpha}^3H(0,0)\gamma_0\intr\frac{|y|^{2(\alpha-1)}}{(1+|y^{\alpha }-b_0|^2)^2}\frac{1-|y^{\alpha }-b_0|^{2}}{1+|y^{\alpha }-b_0|^{2}}dy
\\ &\;\;\;\;
-32{\alpha}^3\log \de_n\intr \frac{|y|^{2(\alpha-1)}}{(1+|y^{\alpha }-\de_n^{-\al}b_n|^2)^2}\tilde\phi_ndy +o(1).\end{aligned}
$$
Then by step 3, \eqref{id1}-\eqref{id2},
\beq\label{ti3}\begin{aligned}\la_n\into V(x) |x|^{2(\alpha-1)}e^{PW_{\la_n}}\phi_nPW_{\la_n} dx&=8\pi{\alpha}\gamma_0+o(1).\end{aligned}
\eeq
Finally, taking into account that $PW_{\la}=O(|\log\la|)$,  we have $\into |\nabla PW_{\la}|^2=\into |x|^{2(\alpha-1)}e^{W_{\la}}PW_{\la}
=O(|\log\la|)$, by which, owing to \eqref{inv2},
\beq\label{ti4} \into |\nabla h_n|\,| \nabla PW_{\la_n}| dx\leq \|h_n\|\, \|PW_{\la_n}\|=o(1).\eeq
By inserting \eqref{ti1}, \eqref{ti3}, \eqref{ti4} into \eqref{ti0} and passing to the limit we deduce $\gamma_0=0$.

\medskip

\noindent{\em Step 5. End of the proof.}

\medskip

We will show that a contradiction arises. According to Step 2 and Step 4 we have $$\tilde\phi_n\to 0\;\ \hbox{  weakly in }\rm{H}_{\al} (\rr^2)\hbox{ and strongly in }\rm{L}_{\al} (\rr^2). $$
By Lemma \ref{aux} 
$$\la_n\into V(x)|x|^{2(\alpha-1)}e^{PW_{\la_n}}\phi_n^2dx=\into |x|^{2(\alpha-1)}e^{W_{\la_n}}\phi_n^2dx+o(1)\leq C\|\tilde\phi_n\|_{\rm L_{\alpha}}^2+o(1)=o(1).$$
Moreover, by \eqref{inv2},
 $$\into \nabla h_n \nabla \phi_ndx=o(1).$$
We multiply the equation in \eqref{lla} by $\phi_n$, we integrate over $\Omega$ and we obtain 
$$\into |\nabla \phi_n|^2dx =\la_n\into V(x) |x|^{2(\alpha-1)}e^{PW_{\la_n}}\phi_n^2 dx -\into \nabla h_n \nabla \phi_ndx =o(1),$$ in contradiction with \eqref{inv2}.

 \end{proof}
 
In addition to \eqref{lla},  let us consider the following linear problem: given
$ h\in {H}_0^1(\Omega)$,  find a function $\phi\in  H^1_0(\Omega)$ and constants $c_1,c_2\in\R$ satisfying
\begin{equation}\label{lla2}
\left\{\begin{aligned}&-\Delta \phi   -\la V(x)|x|^{2(\al-1)} e^{P{W}_\la}\phi=\Delta h+\sum_{j=1,2}c_j Z^j_{\la} |x|^{2(\alpha-1)}e^{W_{\la}}\\ &\into \nabla \phi\nabla PZ^j_{\la}dx=0\;\;j=1,2\end{aligned}\right..
\end{equation}

In order to solve problem \eqref{lla2}, we need to establish an a priori estimate analogous to that of Proposition \ref{inv}. 
\begin{prop}\label{linear2} Let $r>0$ be fixed. 
There exist $\lambda_0>0$ and $C>0$ such that for any $\la \in(0, \la_0)$, any $b\in\R^2$ with $|b|<r\sqrt{\la}$  and   any $h\in  H^1_0(\Omega)$, if 
$(\phi,c_1,c_2)\in  H^1_0(\Omega)\times \R^2$  solves \eqref{lla}, then the following holds $$\|\phi\|    \leq C  |\log\la |    \|h\| .$$ 
\end{prop}
\begin{proof} First observe that 
by \eqref{pzi}
\beq\label{zeta1}\begin{aligned}\into \nabla PZ_\la^1\nabla PZ_\la^2dx&=\into |x|^{2(\alpha-1)}e^{W_{\la}}Z^1_{\la} PZ^2_{\la}dx=\intr |x|^{2(\alpha-1)}e^{W_{\la}}Z^1_{\la} Z^2_{\la}dx+o(1)\\ &=\intr \nabla PZ_1\nabla PZ_2dy +o(1)=o(1).\end{aligned}\eeq Moreover
\beq\label{zeta2}\begin{aligned}\|PZ^1_{\la}\|^2&=\into |x|^{2(\alpha-1)}e^{W_{\la}}Z^1_{\la} PZ^1_{\la}dx=\into |x|^{2(\alpha-1)}e^{W_{\la}}(Z^1_{\la})^2dx+o(1)\\ &=8{\alpha}^2 \intr|y|^{2(\alpha-1)} \frac{|\Re (y^{\alpha }-\de^{-\al}b)|^2}{(1+|y^{\alpha }-\de^{-\al}b|^{2})^4}dy+o(1)=\frac23 \pi{\alpha}+o(1) \end{aligned}\eeq where we have used \eqref{id3}. Similarly 
\beq\label{zeta3}\begin{aligned}\|PZ^2_{\la}\|^2=\frac23 \pi{\alpha}+o(1) .\end{aligned}\eeq
Then, taking into account that $-\Delta PZ^j_{\la} =|x|^{2(\alpha-1)}e^{W_{\la}}Z^j_{\la} $, according to Proposition \ref{inv} we have 
\beq\label{marc}\|\phi\|\leq C\log\la\big(\|h\|+|c_1|+|c_2|\big).\eeq Hence it suffices to estimate the values of the constants $c_j$.
We multiply the  equation in
\eqref{lla2} by $PZ^1_{\la}$ and we find
\beq\label{est}\into \phi  |x|^{2(\alpha-1)}e^{W_{\la}}Z^1_{\la}dx-  \la\into V(x) |x|^{2(\al-1)} e^{P{W}_\la}\phi PZ^1_{\la}dx=\frac23 \pi{\alpha} c_1+o(c_1)+o(c_2)+ O(\|h\|).\eeq
Let us fix $p\in  (1,+\infty)$ sufficiently close to 1. Then, by \eqref{pzi} and \eqref{judo} we may estimate 

$$\begin{aligned}\into  |\phi|  |x|^{2(\alpha-1)}e^{W_{\la}}|PZ^1_{\la}-&Z^1_{\la}|dx\leq C\sqrt\la\into |\phi|  |x|^{2(\alpha-1)}e^{W_{\la}}dx\leq  C\sqrt\la\|\phi\|\, \| |x|^{2(\alpha-1)}e^{W_{\la}}\|_p\\ &\leq C\la^{\frac12 +\frac{1-p}{\al p}}\|\phi\|\leq C\la^{\frac{1}{\al p}}\|\phi\|\end{aligned}$$
and, since $PZ^1_{\la}=O(1)$, using Lemma \ref{aux},
$$\begin{aligned}\into |\phi| \big| |x|^{2(\alpha-1)}&e^{W_{\la}}- \la V(x)|x|^{2(\al-1)} e^{P{W}_\la}\big||PZ^1_{\la}|dx\\&\leq C\into |\phi|  \big||x|^{2(\alpha-1)}e^{W_{\la}}-  \la V(x) |x|^{2(\al-1)} e^{P{W}_\la}\big| dx\leq C\la^{\frac{1}{\al p}}\|\phi\|
.\end{aligned}$$
By inserting the above two estimates into \eqref{est} we obtain 

$$|c_1|+o(c_2)\leq C\|h\|+C\la^{\frac{1}{\al p}}\|\phi\| .$$
We multiply the  equation in
\eqref{lla2} by $PZ^2_{\la}$ and, by a similar argument as above,  we find$$|c_2|+o(c_1)\leq C\|h\|+C\la^{\frac{1}{\al p}}\|\phi\|,$$ and so $$|c_1|+|c_2|\leq C\|h\|+C\la^{\frac{1}{\al p}}\|\phi\|.$$ Combining this with \eqref{marc} we obtain the thesis. 

\end{proof}

\section{The nonlinear problem: a contraction argument}
In order to solve \eqref{eq}, let us consider the following intermediate problem:

\beq\label{inter}\left\{\begin{aligned}&-\Delta(PW_\la+\phi)-\la V(x) |x|^{2(\alpha-1)}e^{PW_\la+\phi}=\sum_{j=1,2}c_j Z_\la^j|x|^{2(\alpha-1)} e^{W_\la},\\ &\phi \in H^1_0(\Omega),\;\;\;\; \into \nabla \phi\nabla PZ_\la^jdx=0 \; \; j=1,2.\end{aligned}\right.\eeq
 
 Then it is convenient to solve as a first step the problem for $\phi$ as a function of $b$. 
 To this aim, first let us rewrite problem \eqref{inter} in a more convenient way. 
 
 For any $ p>1,$ let $$i^*_{p}:L^{p}(\Omega)\to H^1_0(\Omega)$$ be the
adjoint operator of the embedding
$i_{p}:H^1_0(\Omega)\hookrightarrow L^{p\over p-1 }(\Omega),$ i.e.
$u=i^*_{p}(v)$ if and only if $-\Delta u=v$ in $\Omega,$ $u=0$ on
$\partial\Omega.$ We point out that $i^*_{p}$ is a continuous
mapping, namely
\begin{equation}
\label{isp} \|i^*_{p}(v)\| \le c_{p} \|v\|_{p}, \ \hbox{for any} \ v\in L^{p}(\Omega),
\end{equation}
for some constant $c_{p}$ which depends on $\Omega$ and $p.$
Next let us set
$$      {K }:=\hbox{span}\left\{PZ^1_{\la},\ PZ^2_{\la}\right\}$$
and
$$      {K ^\perp }:= \left\{\phi\in H^1_0(\Omega)\ :\ \into \nabla \phi \nabla PZ^1_{\la}dx= \into \nabla \phi \nabla PZ^2_{\la}dx= 0\right\} $$
and  denote by
$$     \Pi : H^1_0(\Omega)\to       {K },\qquad      {\Pi ^\perp}: H^1_0(\Omega)\to       {K ^\perp }$$
the corresponding projections.
Let $      L: K^\perp \to K^\perp$  be the
linear operator defined by
\beq\label{elle}
      L(     \phi):=  \Pi^{\perp}\Big(   {i^*_{p}}\big(       \la V(x)|x|^{2(\alpha-1)}e^{PW_\la}     \phi \big)\Big) - \phi. 
\eeq
Notice that problem \eqref{lla2} reduces to $$L(\phi)=\Pi^\perp h, \quad \phi\in K^\perp.$$
 
 As a consequence of Proposition \ref{linear2} we derive the invertibility of $L$.

\begin{prop}\label{ex} Let $r>0$ be a fixed    number. For any $p>1$ there exist $\lambda_0>0$ and $C>0$ such that for any $\la \in(0, \la_0)$, any $b\in\R^2$ with $\|b\|<r\sqrt{\la}$ and 
 any $h\in K^\perp$  there is a unique solution $ \phi\in K^\perp$ to the problem $$L(\phi)=h.$$ In
particular, $L$ is invertible; moreover, $$\| L^{-1} \| \leq C
|\log \lambda |.$$

\end{prop}
\begin{proof}  Observe that the operator $\phi\mapsto \Pi^\perp\big( {i^*_{p}}\(       \la V(x)|x|^{2(\alpha-1)}e^{PW_\la}     \phi \)\big)$ is a compact operator in $K^\perp$.
   Let us consider the case $h=0$, and take $\phi\in K^\perp$ with $L(\phi)=0$. In other words,  $\phi$ solves
the system \eqref{lla2} with $h=0$ for some $c_1,c_2\in\R$. Proposition \ref{linear2} implies $\phi\equiv 0$. Then, Fredholm's alternative implies the existence and uniqueness result.

Once we have existence, the norm estimate follows directly from Proposition \ref{linear2}.
\end{proof}
Now we come back to our goal of finding a solution to problem \eqref{inter}. In what follows we denote by $N:K^\perp\to K^\perp$ the nonlinear operator 
$$N(\phi)=\Pi^\perp\({i^*_{p}}\big(       \la V(x)|x|^{2(\alpha-1)}e^{PW_\la}(e^{\phi}-1-    \phi) \big)\)$$
 Therefore problem \eqref{inter} turns out to be equivalent to the problem
 
 \beq\label{interop}  L(\phi)+N(\phi)=\tilde R,\quad \phi\in K^\perp\eeq
 where, recalling Lemma \ref{aux},  $$\tilde R=\Pi^\perp\({i^*_{p}}\big(R_\la\big)\)=
 \Pi^\perp\(PW_\la -{i^*_{p}}\big(\la|x|^{2(\alpha-1)}e^{P W_\la}\big)\).$$
 
 We need the following auxiliary lemma.
 
 \begin{lemma}\label{auxnonl} Let $r>0$ be a fixed    number. 
 For any $p> 1$ there exists $\la_0>0$ such that   for any $\la\in (0,\la_0)$, any $b\in \R^2$ with $|b|\leq r \sqrt{\la}$  and any $\phi_1,\phi_2\in H_0^1(\Omega)$ with $\|\phi\|_1,\,\|\phi_2\|<1$ the following holds
\beq\label{skate1}\|e^{\phi_1}-\phi_1-e^{\phi_2}+\phi_2\|_p\leq C(\|\phi_1\|+\|\phi_2\|)\|\phi_1-\phi_2\|,\eeq
 \beq\label{skate2}\|N(\phi_1)-N(\phi_2)\|\leq C\la^{\frac{1-p^2}{{\alpha} p^2}}(\|\phi_1\|+\|\phi_2\|)\|\phi_1-\phi_2\|.\eeq
 
 \end{lemma}
 \begin{proof}
A straightforward computation give that  the inequality $|e^a-a-e^b+b|\leq e^{|a|+|b|}(|a|+|b|)|a-b|$ holds  for all $a,b\in \R$. Then, by applying H\"older's inequality with $\frac1q+\frac1r+\frac1t=1$, we derive
$$\|e^{\phi_1}-\phi_1-e^{\phi_2}+\phi_2\|_p\leq C\|e^{|\phi_1|+|\phi_2|}\|_{pq}(\|\phi_1\|_{pr}+\|\phi_2\|_{pr})\|\phi_1-\phi_2\|_{pt}$$
 and \eqref{skate1} follows by using Lemma \ref{tmt} and the continuity of the embeddings $H^1_0(\Omega)\subset L^{pr}(\Omega)$ and $H^1_0(\Omega)\subset L^{pt}(\Omega)$. 
 Let us prove \eqref{skate2}. According to \eqref{isp}  we get
 $$\|N(\phi_1)-N(\phi_2)\|\leq C\|\la V(x)|x|^{2(\alpha-1)} e^{PW_\la}(e^{\phi_1}-\phi_1-e^{\phi_2}+\phi_2)\|_p,$$
and by H\"older's inequality with $\frac1p+\frac1q=1$ we derive 
$$\begin{aligned}\|N(\phi_1)-N(\phi_2)\|&\leq C\|\la |x|^{2(\alpha-1)}e^{PW_\la}\|_{p^2}\|e^{\phi_1}-\phi_1-e^{\phi_2}+\phi_2|\|_{pq}\\ &\leq C\|\la |x|^{2(\alpha-1)}e^{PW_\la}\|_{p^2}(\|\phi_1\|+\|\phi_2\|)\|\phi_1-\phi_2\|\end{aligned}
$$ by \eqref{skate1}, and  the conclusion follows recalling \eqref{judo}. 
 \end{proof}
Problem \eqref{inter} or, equivalently, problem \eqref{interop}, turns out to be solvable for any choice of point $b$ with $|b|\leq r\sqrt{\la}$, provided
that $\la$ is sufficiently small. Indeed we have the following result.

\begin{prop}\label{nonl} Let $r>0$ be fixed. 
For any $\e\in (0,\frac1\al)$ there exists $\la_0>0$ such that for any $\la\in (0,\la_0)$ and any $b\in \R^2$ with $|b|<r\sqrt\la$ there is a unique $\phi_\la=\phi_{\la,b}\in K^\perp$ satisfying \eqref{inter} for some $c_1,c_2\in \R$ and 
$$\|\phi_{\la}\|\leq \la^{\frac{1}{\alpha}-\e}.$$ 
\end{prop}
\begin{proof} Since, as we have observed, problem \eqref{interop} is equivalent to problem \eqref{inter}, 
we will show that problem \eqref{interop} can be solved via a contraction mapping argument. Indeed, in virtue of Proposition \ref{ex}, let us introduce the map
$$T:=L^{-1}(\tilde R-N(\phi)),\quad \phi\in K^\perp.$$
 Let us fix $$0<\eta<\min\Big\{\e,\frac1\alpha-\e\Big\}$$ and $p>1$ sufficiently close to 1. According to \eqref{isp} and Lemma \ref{aux} we have
\beq\label{non1}\|\tilde R\|=O(\la^{\frac1\alpha-\eta}).\eeq
Similarly, by \eqref{skate2}, choosing $p>1$ sufficiently close to 1, we get 
\beq\label{non2}\|N(\phi_1)-N(\phi_2)\|\leq C\la^{-\eta}(\|\phi_1\|+\|\phi_2\|)\|\phi_1-\phi_2\|\quad \forall \phi_1,\phi_2\in H_0^1(\Omega), \|\phi_1\|,\|\phi_2\|<1.\eeq In particular, by taking $\phi_2=0$, 
\beq\label{non3}\|N(\phi)\|\leq C\la^{-\eta}\|\phi\|^2\quad \forall \phi\in H_0^1(\Omega), \|\phi\|<1.\eeq

We claim that $T$ is a contraction map over the ball $$\Big\{\phi\in K^\perp\,\Big|\, \|\phi\|\leq \la^{\frac1\alpha-\e}\Big\}$$ provided that $\la$ is small enough. Indeed, combining  Proposition \ref{ex}, \eqref{non1}, \eqref{non2}, \eqref{non3} with the choice of $\eta$,
we have 
$$\|T(\phi)\|\leq C|\log\la|(\la^{\frac1\alpha-\eta}+\la^{-\eta}\|\phi\|^2)<\la^{\frac1\alpha-\e},$$
$$\begin{aligned}\|T(\phi_1)-T(\phi_2)\|&\leq C|\log\la|\|N(\phi_1)-N(\phi_2)\|\leq C\la^{-\eta}|\log\la| (\|\phi_1\|+\|\phi_2\|)\|\phi_1-\phi_2\|\\ &<\frac12\|\phi_1-\phi_2\|.\end{aligned}$$

\end{proof}

\section{
Proof of Theorems \ref{th1}-\ref{th2} and Theorem \ref{main1}-\ref{main2}}
After problem \eqref{inter} has been solved according to Proposition \ref{nonl}, then we find a solution to the original problem \eqref{proreg} if $b$   is such that $$c_j=0\hbox{ for }j=1,2.$$ 
 Let us find the condition satisfied by $b$ in order to get the $c_j$'s equal to zero. 

\subsection*{Proof of Theorems \ref{main1}-\ref{main2}} 

We multiply the equation in \eqref{inter} by   $PZ_\la^j$ and integrate over $\Omega$:
\beq\label{masca}\begin{aligned}\into \nabla (PW_\la+\phi_{\la}) \nabla PZ^j_{\la} dx &-\la \into V(x) |x|^{2(\alpha-1)}e^{P W_\la+\phi_{\la}}PZ^j_{\la} dx\\&=\sum_{h=1,2}c_h \into Z_\la^h|x|^{2(\alpha-1)} e^{W_\la}PZ_\la^j dx.\end{aligned}
\eeq
The object  is now to expand each integral of the above identity and analyze the  leading term.
 Let us begin by observing that  the orthogonality in \eqref{inter} gives
\beq\label{masca1} \into \nabla \phi_{\la} \nabla PZ^j_{\la} dx =\into |x|^{2(\alpha-1)} e^{W_\la} \phi_{\la}Z_\la^j dx=0\eeq
and, by \eqref{zeta1}-\eqref{zeta2}, \beq\label{masca2}\into Z_\la^h|x|^{2(\alpha-1)} e^{W_\la}PZ_\la^j dx= \into \nabla PZ_\la^h \nabla PZ^j_{\la} dx=\left\{\begin{aligned}&\frac23 {\pi{\alpha}}+o(1) &\hbox{ if }h=j\\ &o(1)&\hbox{ if }h\neq j\end{aligned}\right..\eeq
Using the expansion \eqref{cuore} 
we get
\beq\label{coll}\begin{aligned}&\into \nabla PW_\la \nabla PZ^j_{\la} dx -\la \into V(x) |x|^{2(\alpha-1)}e^{P W_\la}PZ^j_{\la} dx\\ &=\into |x|^{2(\alpha-1)}  e^{W_\la} PZ^j_{\la} dx-\la \into V(x)|x|^{2(\alpha-1)}e^{P W_\la}PZ^j_{\la} dx\\ &=\into |x|^{2(\alpha-1)}e^{W_\la} \Big(1-\frac{a(x)}{a(0)}e^{-4\pi({\alpha-1}) (H(x,0)-H(0,0))+8\pi\sum_{i=0}^{\al-1}(H(x, \beta_i) -H(0,\beta_i))+O(\de^{2\al})}\Big)PZ^j_{\la} dx
.\end{aligned}\eeq
Recalling  \eqref{expa} and Lemma \ref{robin10}  we deduce 
$$\begin{aligned}&-4\pi({\alpha}-1) (H(x,0)-H(0,0))+8\pi\sum_{i=0}^{\al-1}(H(x, \beta_i)-H(0,\beta_i))\\ &= 4\pi(\al+1) \sum _{k=1}^{2}\frac{1}{k!}\Re\Big(\frac{d^k\tilde H}{dx^k} (0,0)x^k\Big)
+\frac{8\pi}{(\al-1)!}\Re\bigg(\frac{\partial^{\al+1}\tilde H}{\partial p^\al\partial x}(0,0)bx\bigg)\\ &\;\;\;\;+O(|b||x|^2)+O(|b|^2|x|)+O(|x|^3).\end{aligned}
$$
Consequently using the Taylor expansion $e^y=1+y+\frac{y^2}{2}+O(|y|^3)$,  $$\begin{aligned}&e^{-4\pi({\alpha-1}) (H(x,0)-H(0,0))+8\pi\sum_{i=0}^{\al-1}(H(x, \beta_i) -H(0,\beta_i))+O(\de^{2\al})}\\ &=1+4\pi(\al+1) \sum _{k=1}^{2}\frac{1}{k!}\Re\Big(\frac{d^k\tilde H}{dx^k} (0,0)x^k\Big)+\frac{8\pi}{(\al-1)!}\Re\bigg(\frac{\partial^{\al+1}\tilde H}{\partial p^\al\partial x}(0,0)bx\bigg)\\ &\;\;\;\;
+\frac12\bigg(4\pi(\al+1) \sum _{k=1}^{2}\frac{1}{k!}\Re\Big(\frac{d^k\tilde H}{dx^k} (0,0)x^k\Big)\bigg)^2\\ &\;\;\;\;+O(|b||x|^2)+O(|b|^2|x|)+O(|x|^3)+O(\de^{2\al})
\\ &=1
+4\pi(\al+1) \sum _{k=1}^{2}\frac{1}{k!}\Re\Big(\frac{d^k\tilde H}{dx^k} (0,0)x^k\Big)+\frac{8\pi}{(\al-1)!}\Re\bigg(\frac{\partial^{\al+1}\tilde H}{\partial p^\al\partial x}(0,0)bx\bigg)\\ &\;\;\;\;+8\pi^2(\al+1)^2 \bigg(\Re\Big(\frac{d\tilde H}{dx} (0,0)x\Big)\bigg)^2\\ &\;\;\;\; +O(|b||x|^2)+O(|b|^2|x|)+O(|x|^3)+O(\de^{2\al}) .
\end{aligned}$$ 
By assumptions \eqref{assurd0} and \eqref{assurd}-\eqref{assurd2}  in Theorems \ref{th1} and \ref{th2}, respectively, taking into account that $\Re(\frac{d\tilde H}{dx} (0,0)x)=\langle \nabla_xH(0,0), x\rangle$,    we get $$\begin{aligned}\frac{a(x)}{a(0)}&=1+\frac{\langle \nabla a(0), x\rangle}{a(0)}
+\frac{1}{2a(0)}\Big( a_{11}x_1^2+a_{22}x_2^2
\Big)+O(|x|^3)\\ &=1-4\pi(\al+1)\Re\Big(\frac{d\tilde H}{dx} (0,0)x\Big)
+\frac{1}{2a(0)}\Big( a_{11}(\Re \,x)^2+a_{22}(\Im \,x)^2
\Big)+O(|x|^3),
\end{aligned}$$ 
and then we derive
 \beq\label{insert}\begin{aligned}&\frac{a(x)}{a(0)}e^{-4\pi({\alpha-1}) (H(x,0)-H(0,0))+8\pi\sum_{i=0}^{\al-1}(H(x, \beta_i) -H(0,\beta_i))+O(\de^{2\al})}\\ &=1
+2\pi(\al+1)\Re\Big(\frac{d^2\tilde H}{dx^2} (0,0)x^2\Big)
+\frac{8\pi}{(\al-1)!}\Re\bigg(\frac{\partial^{\al+1}\tilde H}{\partial p^\al\partial x}(0,0)bx\bigg)\\ &\;\;\;\;
-8\pi^2(\al+1)^2 \bigg(\Re\Big(\frac{d\tilde H}{dx} (0,0)x\Big)\bigg)^2\\ &
 \;\;\;\; +\frac{1}{2a(0)}\Big( a_{11}(\Re \,x)^2+a_{22}(\Im \,x)^2
 \Big)+O(|b||x|^2)+O(|b|^2|x|)+O(|x|^3)+O(\de^{2\al}). \end{aligned}\eeq
First let us assume that $\al\geq 3$: let us insert the above expansion into \eqref{coll} and, using Lemma \ref{copy0cor} and next Corollary \ref{corocoro} we get

\beq\label{coll2}\begin{aligned}&\into \nabla PW_\la \nabla PZ^j_{\la} dx -\la \into V(x) |x|^{2(\alpha-1)}e^{P W_\la}PZ^j_{\la} dx\\ &=8\pi^2(\al+1)^2\into |x|^{2(\alpha-1)}e^{W_\la}PZ_\la^j 
\bigg(\Re\Big(\frac{d\tilde H}{dx} (0,0)x\Big)\bigg)^2dx\\ &\;\;\;\;
-\frac{1}{2a(0)}\into |x|^{2(\alpha-1)}e^{W_\la}PZ_\la^j \Big( a_{11}(\Re \,x)^2+a_{22}(\Im \,x)^2\Big)dx \\ &\;\;\;\;+ O(\de^3)+O(\de^2|b|)+O(\de|b|^2)\\ &=
\frac{\de^2}{2} \Bigg(8\pi^2(\al+1)^2\Big|\frac{d\tilde H}{dx} (0,0)\Big|^2- \frac{a_{11}+a_{22}}{2a(0)}\bigg) \intr |x|^{2\al}e^{W_\la}Z^j_\la dx\\ &\;\;\;\;+ O(\de^3)+O(\de^2|b|)+O(\de|b|^2).\end{aligned}\eeq

We have thus obtained that if $\al\geq 3$ then 
\beq\label{masca3}\begin{aligned} &\into \nabla PW_\la \nabla PZ^j_{\la} dx  -\la \into V(x) |x|^{2(\alpha-1)}e^{P W_\la}PZ^j_{\la} dx\\ &=A\de^2
F_j(\de^{-\al}b) 
+ O(\de^3)+O(\de^2|b|)+O(\de|b|^2)
\end{aligned}\eeq
where $$A:= 4\pi^2(\al+1)^2\Big|\frac{d\tilde H}{dx} (0,0)\Big|^2- \frac{a_{11}+a_{22}}{4a(0)}\neq 0$$ thanks to assumptions \eqref{assurd0} 
in Theorem \ref{th1} and $F=(F_1, F_2)$ is the map defined in Lemma \ref{finalaux}.

Next assume that $\al=2$.  If $\Omega$ is $\ell$-symmetric for some $\ell\geq3$ in the sense of \eqref{assurd}, then $\tilde H(x,0)$ is $3$-symmetric too:
$$\tilde H(e^{{\rm i}\frac{2\pi}{\ell}}x,0)= \tilde H(x,0)\quad \forall x\in\Omega;$$
 this implies that its Taylor expansion at $0$ involves only the powers corresponding to integers multiples of $\ell$ and, consequently,  
 $$\frac{d\tilde H}{dx} (0,0)=\frac{d^2\tilde H}{dx^2} (0,0)=0.$$
 Then let us insert \eqref{insert} into \eqref{coll} and, using Lemma \ref{copy0cor} and next Corollary \ref{corocorocoro} we get for $j=1$
$$\begin{aligned}&\into \nabla PW_\la \nabla PZ^1_{\la} dx -\la \into V(x) |x|^{2(\alpha-1)}e^{P W_\la}PZ^1_{\la} dx\\ &=
-\frac{1}{2a(0)}\into |x|^{2(\alpha-1)}e^{W_\la}PZ_\la^1 \Big( a_{11}(\Re \,x)^2+a_{22}(\Im \,x)^2
\Big)dx + O(\de^3)+O(\de^2|b|)+O(\de|b|^2)\\ &=
- \de^2\frac{a_{11}+a_{22}}{4a(0)}\intr |x|^{2\al} e^{W_\la}Z^1_\la dx-\pi\al^2\de^2\frac{a_{11}-a_{22}}{a(0)}
 + O(\de^3)+O(\de^2|b|)+O(\de|b|^2)\end{aligned}$$
and, similarly for $j=2$
$$\begin{aligned}&\into \nabla PW_\la \nabla PZ^2_{\la} dx -\la \into V(x) |x|^{2(\alpha-1)}e^{P W_\la}PZ^2_{\la} dx\\ &=
-\frac{1}{2a(0)}\into |x|^{2(\alpha-1)}e^{W_\la}PZ_\la^2 \Big( a_{11}(\Re \,x)^2+a_{22}(\Im \,x)^2
\Big)dx + O(\de^3)+O(\de^2|b|)+O(\de|b|^2)\\ &=
- \de^2\frac{a_{11}+a_{22}}{4a(0)} \intr |x|^{2\al} e^{W_\la}Z^2_\la dx
 + O(\de^3)+O(\de^2|b|)+O(\de|b|^2).\end{aligned}$$
 Therefore, using \eqref{assurd2} we conclude that \eqref{masca3} holds for any $\al\geq 2$ for some $A\neq 0$.

Finally let us fix $\e>0$ sufficiently small and $p>1$ sufficiently close to 1. 
Next let  $1<q<\infty$ be such that  $\frac1p+\frac1q=1$. Then, recalling that $\de^{2\al}\sim\la$ according to \eqref{delta},
\eqref{skate1} with $\phi_2=0$ and Proposition \ref{nonl} give
$$\|e^{\phi_\la}-1-\phi_\la\|_q\leq  C\|\phi\|^2\leq \de^{4-4\alpha \e}  $$
and, consequently, 
\beq\label{tmtsur}\|e^{\phi_\la}-1\|_q\leq C\|\phi_\la\|\leq \de^{2-{2\alpha} \e}  .\eeq
Therefore,  the orthogonality  \eqref{masca1}  and Lemma \ref{aux} imply 
 $$\begin{aligned} 
  \into |x|^{2(\alpha-1)}e^{W_\la}(e^{\phi_\la}-1)Z^j_{\la} dx
&  =\into |x|^{2(\alpha-1)}e^{W_\la}(e^{\phi_\la}-1-\phi_\la)Z^j_{\la} dx\\ &
= O(\|  |x|^{2(\alpha-1)}e^{ W_\la}(e^{\phi_\la}-1-\phi_\la)\|_1)\\ &=
O(\|   e^{W_\la} |x|^{2(\alpha-1)}\|_p\|e^{\phi_\la}-1-\phi_\la\|_q)\\ &=O(\de^{\frac{2}{p}-2}  \de^{4-4\alpha\e})
\end{aligned}$$ 
and,  by using again  Lemma \ref{aux} and \eqref{pzi},
\beq\label{masca4}\begin{aligned} \la \into V(x) |x|^{2(\alpha-1)}e^{P W_\la}(e^{\phi}-1)PZ^j_{\la} dx&=\into |x|^{2(\alpha-1)}e^{W_\la}(e^{\phi}-1)Z^j_{\la} dx+ O(\de^{\frac{2}{p}+2-{2\alpha}\e})
\\ &=O(\de^{\frac{2}{p}-2}  \de^{4-4\alpha \e})+ O(\de^{\frac{2}{p}+2-{2\alpha}\e})=o(\de^3)
\end{aligned}\eeq
provided that  $\e$ is chosen   sufficiently close to  0 and $p$ sufficiently close to 1.

In order to conclude, combining \eqref{masca1}, \eqref{masca2}, \eqref{masca3},  \eqref{masca4},  
 the identities \eqref{masca} turn
out to be equivalent to the system
\beq\label{tccone}\begin{aligned} &\de^2A\big(F_1(\de^{-\al}b)+O(\de)
\big)=\frac{2}{3}\pi\al c_1+o(c_1)+o(c_2), \\ & \de^2A\big( F_2(\de^{-\al}b)+O(\de)
\big)=\frac{2}{3}\pi\al c_2+o(c_1)+o(c_2)\end{aligned}\eeq uniformly for $|b|\leq \de^\al$. 
According to Lemma \ref{finalaux}  we have  $F(0,0)=(0,0)$ and $\det F(0,0)\neq0.$ Then the local invertibility theorem assures that $F$ is invertible in a small ball $B_r$ with center $0$ or, equivalently,   $F(\de^{-\al}b)$ is  invertible in a the  ball $B_{r\de^{\al}}$,   and hence $\deg(F(\de^{-\al}b), B_{r\de^{\al}},0)=1.$ 
Taking into account that $|F(\de^{-\al}b)|\geq c$ for $|b|=r\de^\al$,  the continuity property  of the topological
degree gives that $\deg(F(\de^{-\al}b)+O(\de) , B_{r\de^\al},0)>0$ for $\de$ (hence $\la$) small enough. Then for such $\de$ there exists $b\in B_{r}$ such that 
$$F(\de^{-\al}b)+O(\de)=0.$$
and so the linear system \eqref{tccone}
 has only the trivial solution $c_1=c_2=0$.  Finally  $\de^{-\al}|b|\sim |F(\de^{-\al}b)|=O( \delta) ,$ hence $|b|=O(\de^{\al+1})$.  That concludes the
proof of Theorems \ref{main1}-\ref{main2}.

\bigskip

\noindent{\bf{Proof of Theorems \ref{th1}-\ref{th2}.} }Theorems \ref{main1}-\ref{main2} provide a solution to the  problem \eqref{proreg}  of the form $$v_\la=PW_\la+\phi_\la$$   for some $b=b_\la$ with $|b_\la|=O(\de^{\al+1})$. 
So we have \beq\label{exan0}\sum_{k=0}^{\al-1} H(x,{\beta_i}_\la)=\al H(x,0)+O(\de),\eeq \beq\label{exan}\log (\de^{2\al}+|x^{\al}+ b_\la|^2)
=\log(\de^{2\al}+|x|^{2\al})+O(\de)\eeq
uniformly for $x\in\overline\Omega$.
Clearly, by \eqref{chva},  $$u_\la=v_\la-4\pi (\al-1) G(x,0)$$ solves equation \eqref{eq} and \eqref{the1} of Theorem \ref{th1} follows from \eqref{pro-exp1} and \eqref{exan0}\eqref{exan}.
Moreover, using \eqref{judo} and \eqref{tmtsur}, by H\"older's inequality with $\frac1p+\frac1q=1$ we get
$$\begin{aligned}\la\||x|^{2(\al-1)}V(x)(e^{{v}_\la}- e^{{PW}_\la})\|_1&=\la\||x|^{2(\al-1)}V(x)e^{{PW}_\la}(e^{{\phi}_\la}-1)\|_{1}
\\ &\le\la \||x|^{2(\al-1)}V(x)e^{{PW}_\la}\|_p\|e^{\phi_\la}-1\|_q\\ &=O(\la^{\frac{1-p}{\al p}+\frac{1}{\al}-\e})=o(1), \end{aligned}$$ if $p$ is chosen sufficiently close to 1 and $\e$ sufficiently close to $0$.
Then, by \eqref{quantum} and Lemma \ref{aux} 
$$\begin{aligned}\la \into a(x)e^{u_\la}dx&=\la \into |x|^{2(\al-1)}V(x)e^{v_\la }dx=\la\into|x|^{2(\al-1)}V(x)e^{{PW}_\la}dx +o(1)\\ &=\into|x|^{2(\al-1)}e^{{W}_\la}dx+o(1)=\intr|x|^{2(\al-1)}e^{{W}_\la}dx+o(1)=8\pi\al+o(1). \end{aligned}$$
Similarly for every neighborhood $U$ of $0$ $$\la \int_U a(x)e^{u_\la}dx\to 8\pi\al. $$
So \eqref{the3} is verified and Theorem \ref{th1} is thus completely proved. 

\begin{lemma}\label{copy0cor} Let $\alpha\geq 2$ and $\xi\in \C$. For any $\gamma=0,1,\ldots, \al-1$ the following holds:
$$\into |x|^{2(\al-1)}e^{W_\la} PZ^j_\la\Re(\xi x^\gamma)dx =O(\de^{\al+\gamma}),\quad \into |x|^{2(\al-1)}e^{W_\la} PZ_\la^j\Im(\xi x^\gamma)dx=O(\de^{\al+\gamma})$$ and 
$$\into |x|^{2(\al-1)}e^{W_\la} PZ^1_\la\Re(\xi x^\al)dx=4\pi\al^2\de^\al\Re(\xi)+O(\de^{2\al})$$ $$\into |x|^{2(\al-1)}e^{W_\la} PZ_\la^1\Im(\xi x^\al)dx=4\pi\al^2\de^\al\Im(\xi)+O(\de^{2\al})$$
$$\into |x|^{2(\al-1)}e^{W_\la} PZ^2_\la\Re(\xi x^\al)dx=-4\pi\al^2\de^\al\Im(\xi)+O(\de^{2\al})$$ $$\into |x|^{2(\al-1)}e^{W_\la} PZ_\la^2\Im(\xi x^\al)dx=4\pi\al^2\de^\al\Re(\xi)+O(\de^{2\al})$$
uniformly for $b$ in a small neighborhood of $0$.
\end{lemma}
\begin{proof} Let us first show the  identities for $j=1$ and $\xi=1$.
By \eqref{pzi} for $\gamma=0,1,\ldots, \al$  we compute
$$\begin{aligned}&\into |x|^{2(\al-1)}e^{W_\la} PZ^1_\la\Re(x^\gamma)dx\\&=8\al^2\de^\gamma\int_{\frac{\Omega}{\de}}\frac{|y|^{2(\al-1)}}{(1+|y^\al-\de^{-\al}b|^2)^3}\Re(y^\al-\de^{-\al}b)\Re(y^\gamma) dy +O(\de^{\al+\gamma}) \\ &=8\al^2\de^\gamma\intr\frac{|y|^{2(\al-1)}}{(1+|y^\al-\de^{-\al}b|^2)^3}\Re(y^\al-\de^{-\al}b)\Re(y^\gamma) dy+O(\de^{\al+\gamma}).\end{aligned} $$
If $\gamma=1,\ldots, \al-1$ the thesis follows from Lemma \ref{copy0}. If $\gamma=0$, then by applying Lemma \ref{copy}
$$\begin{aligned}\intr\frac{|y|^{2(\al-1)}}{(1+|y^\al-\de^{-\al}b|^2)^3}\Re(y^\al-\de^{-\al}b) dy&=\frac{1}{\al}\intr\frac{1}{(1+|y-\de^{-\al}b|^2)^3}\Re(y-\de^{-\al}b) dy \\ &= \frac{1}{\al}\intr\frac{y_1}{(1+|y|^2)^3}dy=0\end{aligned}$$
and we get the first estimate for $\xi=1$. The second estimate with $\xi=1$ is analogous. 
Next, if $\gamma=\alpha$ then again by Lemma \ref{copy} 
$$\begin{aligned}&\intr\frac{|y|^{2(\al-1)}}{(1+|y^\al-\de^{-\al}b|^2)^3}\Re(y^\al-\de^{-\al}b) \Re(y^\al)dy\\&=\frac{1}{\al}\intr \frac{1}{(1+|y-\de^{-\al}b|^2)^3}\Re(y-\de^{-\al}b) \Re(y)dy\\ &
=\frac{1}{\al}\intr \frac{1}{(1+|y|^2)^3}\Re(y) \Re(y+\de^{-\al}b)dy\\ &
= \frac{1}{\al}\intr \frac{1}{(1+|y|^2)^3}y_1(y_1-\de^{-\al}\Re(b))dy\\ &=\frac{1}{\al}\intr \frac{(y_1)^2}{(1+|y|^2)^3} dy -\de^{-\al}\Re(b) \frac{1}{\al}\intr \frac{y_1}{(1+|y|^2)^3} dy
\\ &=\frac{\pi}{2} \end{aligned} $$
since $\intr \frac{(y_1)^2}{(1+|y|^2)^3} dy=\frac{1}{2} \intr \frac{|y|^2}{(1+|y|^2)^3} dy=\frac{\pi}{2} $ and $\intr \frac{y_1}{(1+|y|^2)^3} dy=0$. Similarly $$\begin{aligned}&\intr\frac{|y|^{2(\al-1)}}{(1+|y^\al-\de^{-\al}b|^2)^3}\Re(y^\al-\de^{-\al}b) \Im(y^\al)dy\\&
=\frac{1}{\al}\intr \frac{1}{(1+|y|^2)^3}\Re(y) \Im(y+\de^{-\al}b)dy\\ &
= \frac{1}{\al}\intr \frac{1}{(1+|y|^2)^3}y_1 (y_2+\de^{-\al}\Re(b))dy\\ &=0
.\end{aligned}$$ Taking into account that $$\Re(\xi x^\gamma)=\Re(\xi) \Re(x^\gamma)-\Im(\xi)\Im(x^\gamma),\quad \Im(\xi x^\gamma)=\Re(\xi)\Im(x^\gamma)+\Im(\xi)\Re(x^\gamma)$$ we obtain the thesis for $j=1$ and any $\xi\in \C$. The remaining estimates with $j=2$ are analogous.
\end{proof}

\begin{cor}\label{corocoro} Let $\alpha\geq 3$ and $\xi_1, \, \xi_2\in\C$. Then   $$ \into |x|^{2(\al-1)}e^{W_\la} PZ_\la^j\Re(\xi_1x)\Re(\xi_2 x)dx=\frac{\de^2}{2} \langle \xi_1, \xi_2\rangle\intr |x|^{2}e^{W_\la}Z_\la^jdx
+O(\de^{\al+2}),$$
$$ \into |x|^{2(\al-1)}e^{W_\la} PZ_\la^j\Im(\xi_1x)\Im(\xi_2 x)dx=\frac{\de^2}{2} \langle \xi_1, \xi_2\rangle\intr |x|^{2\al}e^{W_\la}Z_\la^jdx+O(\de^{\al+2})$$
uniformly for $b$ in a small neighborhood of $0$.
\end{cor}
\begin{proof}
Since $(\Re(x))^2= \frac{|x|^2}{2}+\frac{\Re(x^2)}{2}$, according to Lemma \ref{copy0cor} we have 
$$\begin{aligned}&\into |x|^{2(\al-1)}e^{W_\la} PZ_\la^1\big(\Re(x)\big)^2dx\\&=\frac12\into |x|^{2(\al-1)}e^{W_\la} PZ_\la^1|x|^2dx+O(\de^{\al+2})\\ &=4\al^2\de^2\int_{\frac{\Omega}{\de}}\frac{|y|^{2\al}}{(1+|y^\al-\de^{-\al}b|^2)^3}\Re(y^\al-\de^{-\al}b) dy+O(\de^{\al+2}) \\ &=4\al^2\de^2\intr\frac{|y|^{2\al}}{(1+|y^\al-\de^{-\al}b|^2)^3}\Re(y^\al-\de^{-\al}b) dy +O(\de^{\al+2})\\ &
=\frac{\de^2}{2} \intr |x|^{2\al}e^{W_\la}Z_\la^1dx+O(\de^{\al+2})
.
\end{aligned}$$ 
Similarly, using now  $(\Im(x))^2= \frac{|x|^2}{2}-\frac{\Re(x^2)}{2}$,
$$\begin{aligned}\into |x|^{2(\al-1)}e^{W_\la} PZ_\la^1\big(\Im(x)\big)^2dx&=\frac{\de^2}{2} \intr |x|^{2\al}e^{W_\la}Z_\la^1dx +O(\de^{\al+2}).
\end{aligned}$$ 
Moreover, since $\Re(x)\Im(x)=\frac{\Im(x^2)}{2}$, by Lemma \ref{copy0cor}
$$\begin{aligned}\into |x|^{2(\al-1)}e^{W_\la} PZ_\la^1\Re(x)\Im(x)dx= O(\de^{\al+2}).\end{aligned}$$ 
The thesis follows for $j=1$  since $\Re(\xi x)=\Re(\xi) \Re(x)-\Im(\xi)\Im(x)$.  The proof for $j=2$ follows analogously. 

\end{proof}
\begin{cor}\label{corocorocoro} Let $\alpha=2$ and $\xi_1, \, \xi_2\in\C$. Then $$ \begin{aligned}\into |x|^{2(\al-1)}e^{W_\la} PZ_\la^1\Re(\xi_1x)\Re(\xi_2 x)dx&=\frac{\de^2}{2}\langle \xi_1, \xi_2\rangle \int |x|^{2\al} e^{W_\la} Z_\la^1dx\\ &\;\;\;\;+2\pi\al^2\de^2\Re(\xi_1 \xi_2) +O(\de^{\al+2}),\end{aligned}$$
$$ \begin{aligned}\into |x|^{2(\al-1)}e^{W_\la} PZ_\la^2\Re(\xi_1x)\Re(\xi_2 x)dx&=\frac{\de^2}{2}\langle \xi_1, \xi_2\rangle \int |x|^{2\al} e^{W_\la} Z_\la^2dx\\ &\;\;\;\;-2\pi\al^2\de^2\Im(\xi_1 \xi_2) +O(\de^{\al+2}),\end{aligned}$$
$$\begin{aligned} \into |x|^{2(\al-1)}e^{W_\la} PZ_\la^1\Im(\xi_1x)\Im(\xi_2 x)dx&=\frac{\de^2}{2}\langle \xi_1, \xi_2\rangle \int |x|^{2\al} e^{W_\la} Z_\la^jdx \\ &\;\;\;\;-2\pi\al^2\de^2\Re (\xi_1, \xi_2)+O(\de^{\al+2})\end{aligned}$$
$$\begin{aligned} \into |x|^{2(\al-1)}e^{W_\la} PZ_\la^2\Im(\xi_1x)\Im(\xi_2 x)dx&=\frac{\de^2}{2}\langle \xi_1, \xi_2\rangle \int |x|^{2\al} e^{W_\la} Z_\la^jdx \\ &\;\;\;\;+2\pi\al^2\de^2\Im(\xi_1\xi_2)+O(\de^{\al+2})\end{aligned}$$
uniformly for $b$ in a small neighborhood of $0$. 
\end{cor}
\begin{proof}
Since $(\Re(x))^2= \frac{|x|^2}{2}+\frac{\Re(x^2)}{2}$, according to Lemma \ref{copy0cor} we have 
$$\begin{aligned}&\into |x|^{2(\al-1)}e^{W_\la} PZ_\la^1\big(\Re(x)\big)^2dx\\&=\frac12\into |x|^{2(\al-1)}e^{W_\la} PZ_\la^1|x|^2dx+2\pi\al^2\de^2+O(\de^{\al+2})\\ &=4\al^2\de^2\int_{\frac{\Omega}{\de}}\frac{|y|^{2\al}}{(1+|y^\al-\de^{-\al}b|^2)^3}\Re(y^\al-\de^{-\al}b) dy+2\pi\al^2\de^2+O(\de^{\al+2}) \\ &=4\al^2\de^2\intr\frac{|y|^{2\al}}{(1+|y^\al-\de^{-\al}b|^2)^3}\Re(y^\al-\de^{-\al}b) dy +2\pi\al^2\de^2+O(\de^{\al+2})
\\ &=\frac{\de^2}{2}\int |x|^{2\al} e^{W_\la} Z_\la^1dx+2\pi\al^2\de^2+O(\de^{\al+2})
.
\end{aligned}$$ 
Similarly, using now  $(\Im(x))^2= \frac{|x|^2}{2}-\frac{\Re(x^2)}{2}$,
$$\begin{aligned}\into |x|^{2}e^{W_\la} PZ_\la^1\big(\Im(x)\big)^2dx&=\frac{\de^2}{2}\int |x|^{2\al} e^{W_\la} Z_\la^jdx-2\pi\al^2\de^2+O(\de^{\al+2}).
\end{aligned}$$ 
Moreover, since $\Re(x)\Im(x)=\frac{\Im(x^2)}{2}$, by Lemma \ref{copy0cor}
$$\begin{aligned}\into |x|^{2}e^{W_\la} PZ_\la^1\Re(x)\Im(x)dx=
O(\de^{\al+2}).\end{aligned}$$ 
The first estimate follows  since $\Re(\xi x)=\Re(\xi) \Re(x)-\Im(\xi)\Im(x)$ and $\Im(\xi x)=\Re(\xi)\Im(x)+\Im(\xi)\Re(x).$ The remaining estimates  follow analogously.

\end{proof}

\begin{lemma} 
\label{finalaux}Let ${\alpha}\geq 2$ be an integer and let $F:\R^2\to \R^2$ be defined  by$$F(B)=\left(\begin{aligned} &\intr  \frac{|y|^{2\al}}{(1+|y^{\alpha }-B|^2)^3}\Re (y^{\alpha }-B)dy\\ &\intr  \frac{|y|^{2\alpha}}{(1+|y^{\alpha }-B|^2)^3}\Im (y^{\alpha }-B)dy\end{aligned}\right).$$
Then $F(0)=0$ and $det (DF(0))\neq 0$.
\end{lemma}

\begin{proof} According to Lemma \ref{copy} we have 
$$\intr  \frac{|y|^{2\al}}{(1+|y|^{{2\alpha}})^3}\Re (y^{\alpha })dy=\frac{1}{\al}\intr \frac{|y|^{2/\al}}{(1+|y|^2)^3}\Re (y)dy=\frac{1}{\al}\intr \frac{|y|^{2/\al}}{(1+|y|^2)^3}y_1dy=0
.$$ Similarly 
$$\intr  \frac{|y|^{2\al}}{(1+|y|^{{2\alpha}})^3}\Im (y^{\alpha })dy=\frac{1}{\al}\intr \frac{|y|^{2/\al}}{(1+|y|^2)^3}y_2dy=0$$ 
by which we immediately get $F(0,0)=0$.
Moreover $$\begin{aligned}\frac{\partial F_1}{\partial B_2}(0,0)&=6\intr\frac{|y|^{2\al}}{(1+|y^{\alpha }|^2)^4}\Im(y^\al)\Re (y^{\alpha })dy=\frac{6}{\al}\intr \frac{|y|^{2/\al}}{(1+|y|^2)^4}\Im (y)\Re(y) dy 
\\ &=\frac{6}{\al}\intr \frac{|y|^{2/\al}}{(1+|y|^2)^4}y_1y_2 dy=0.\end{aligned}$$ Similarly $\frac{\partial F_2}{\partial B_1}(0,0)=0$.
 So $DF$ is a diagonal matrix.  
We compute
$$\begin{aligned}\frac{\partial F_1}{\partial B_1}(0,0)&=-\intr|y|^{2\al}\frac{1+|y|^{2\alpha}-6(\Re (y^{\alpha }))^2}{(1+|y|^{2\alpha})^4} dy
\\ &
=-\frac{1}{\al}\intr|y|^{\frac{2}{\al}}\frac{1+|y|^{2}-6y_1^2}{(1+|y|^{2})^4} dy
.\end{aligned}$$
Using that $\intr|y|^{\frac{2}{\al}}\frac{y_1^2}{(1+|y|^{2})^4} =\frac12 \intr|y|^{\frac{2}{\al}}\frac{|y|^2}{(1+|y|^{2})^4} $ we get
$$\begin{aligned}\frac{\partial F_1}{\partial B_1}(0,0)&=\frac1\al\intr|y|^{\frac2\al}\frac{2|y|^2-1}{(1+|y|^{2})^4} dy.
\end{aligned}$$

Proceeding similarly as above we get
$$\frac{\partial F_2}{\partial B_2}(0,0)=\frac{\partial F_1}{\partial B_1}(0,0)=\frac1\al\intr|y|^{\frac2\al}\frac{2|y|^2-1}{(1+|y|^{2})^4} dy.$$
Then the thesis will follow once we have proved the nonvanishing of the above integral:
\beq\label{nonv} \intr|y|^{\frac2\al}\frac{2|y|^2-1}{(1+|y|^{2})^4} dy\neq 0.\eeq
To see this, let us first compute
$$\begin{aligned}\frac{1}{2\pi}\intr\frac{2|y|^2-1}{(1+|y|^{2})^4} dy&=\int_0^{+\infty} \rho\frac{2\rho^2-1}{(1+\rho^{2})^4} d\rho\\ &=  2\int_0^{+\infty} \rho\frac{1}{(1+\rho^{2})^3} d\rho-3\int_0^{+\infty} \rho\frac{1}{(1+\rho^{2})^4} d\rho=0
 \end{aligned}
$$ by direct integration. So, using that $(\sqrt2 |y|)^{\frac2\al}\leq 1$ if $2|y|^2-1\leq 0$ and  $(\sqrt2 |y|)^{\frac2\al}> 1$ if $2|y|^2-1> 0$ and 
$$\begin{aligned}\intr|y|^{\frac2\al}\frac{2|y|^2-1}{(1+|y|^{2})^4} dy
= (\sqrt2)^{-\frac2\al}\intr(\sqrt2|y|)^{\frac2\al}\frac{2|y|^2-1}{(1+|y|^{2})^4} dy
>  (\sqrt2)^{-\frac2\al}  \intr\frac{2|y|^2-1}{(1+|y|^{2})^4} dy=0\end{aligned}
$$ and \eqref{nonv} follows.
\end{proof}3

\appendix
\renewcommand{\theequation}{\Alph{section}.\arabic{equation}}
\section{}

This appendix is devoted to deduce some integral identities associated to the change of variables: $x\mapsto x^\al$.

\begin{lemma}\label{copy} 
For any $f\in L^1(\R^2)$ we have that $|y|^{2(\al-1)}f(y^\al)\in L^1 (\R^2) $ and \beq\label{normeqq}\intr |y|^{2(\alpha-1)} f(y^\al) dy=\frac1\al\intr f(y )dy.\eeq
\end{lemma}
\begin{proof} It is sufficient to prove the thesis for a smooth function $f$. Using the  polar coordinates $(\rho,\theta)$ and then applying the change of variables $(\rho',\theta')=(\rho^\al,\al\theta)$
we get
$$\begin{aligned} \intr |y|^{2(\alpha-1)} f(y^\al) dy&=\int_0^{+\infty}d\rho\int_0^{2\pi}  \rho^{2\al-1}f(\rho^\al e^{{\rm i}\al \theta}) d\theta 
\\ &= \frac{1}{\al^2}\int_0^{+\infty}d\rho'\int_0^{2\al \pi}  \rho'f(\rho' e^{{\rm i}\theta'}) d\theta' \\ &=\frac{1}{\al}\int_0^{+\infty}d\rho'\int_0^{2\pi}  \rho'f(\rho' e^{{\rm i}\theta'}) d\theta' =\frac{1}{\al}\intr f(y) dy.\end{aligned}$$
\end{proof}

\begin{lemma}\label{copy0}
Let $\gamma=1, \ldots, \al-1$  and let $f$ be such that $f(y)|y|^{\frac\gamma\alpha}\in L^1(\R^2)$. Then 
$$ \intr  
|y|^{2(\al-1)}f(y^\alpha) \Re (y^{\gamma})dy=\intr  |y|^{2(\al-1)} f(y)
\Im (y^{\gamma})dy=0.$$
\end{lemma}
\begin{proof} Observe first that according to Lemma \ref{copy} we have $|y|^{2(\al-1)}f(y^\alpha) \Re (y^{\gamma})\in L^1(\R^2)$. 
Suppose that $f$ is a smooth function.  Using the polar coordinates $(\rho,\theta)$  we get $$\begin{aligned} \intr  |y|^{2(\alpha-1)}f(y)\Re (y^{\gamma})dy&=\int_0^{+\infty}|\rho|^{2\alpha-1+\gamma}d\rho\int_0^{2\pi}\cos(\gamma\theta)f(\rho^\al e^{{\rm i}\al \theta})d\theta.\end{aligned}$$
On the other hand
$$\begin{aligned}\int_0^{2\pi}\cos(\gamma\theta)f(\rho^\al e^{{\rm i}\al \theta})d\theta &= 
\sum_{k=0}^{\al-1}\int_{\frac{2\pi}{\al}k}^{\frac{2\pi}{\al}(k+1)}\cos(\gamma\theta)f(\rho^\al e^{{\rm i}\al \theta})d\theta\\ &=
\sum_{k=0}^{\al-1}\int_0^{\frac{2\pi}{\al}}\cos\Big(\gamma\Big(\theta+\frac{2\pi}{\al}k\Big)\Big)f(\rho^\al e^{{\rm i}\al \theta})d\theta \\ &=\sum_{k=0}^{\al-1}\cos\Big(\gamma\frac{2\pi}{\al}k\Big)\int_0^{\frac{2\pi}{\al}}\cos(\gamma\theta)f(\rho^\al e^{{\rm i}\al \theta})d\theta\\ &\;\;\;-\sum_{k=0}^{\al-1}\sin\Big(\gamma\frac{2\pi}{\al}k\Big)\int_0^{\frac{2\pi}{\al}}\sin(\gamma\theta)f(\rho^\al e^{{\rm i}\al \theta})d\theta.
\end{aligned}$$ The well known identity  $\sum_{k=0}^{\al-1}e^{{\rm i}\frac{2\pi}{\al}\gamma k}=0$ for all $\gamma=1,\ldots,\al-1$ implies  $$\sum_{k=0}^{\al-1}\cos\Big(\gamma\frac{2\pi}{\al}k\Big)=\sum_{k=0}^{\al-1}\sin\Big(\gamma\frac{2\pi}{\al}k\Big)=0.$$
and the first identity follows. The second identity is analogous..

\end{proof}

\section{}
 In this appendix we carry out some asymptotic expansions involving the regular part $H(x,y)$ of the Green's function.  Recalling that for any fixed $p\in\Omega$ the function $H_p:x\mapsto H(x,p)$ is harmonic in $\Omega$, then it admits a holomorphic extension $$\tilde H_p(x)=H_p(x)+{\rm i}h_p(x)\;\hbox{ in } U,\;\; h_p(p)=0\qquad x\approx x_1+{\rm i} x_2,
$$ where $U$ is   any fixed round neighborhood of $0$. 
Setting $$\tilde H(x,p)=\tilde H_p(x)\quad \forall x,p\in U$$ 
by the symmetry  $H(x,p)=H(p,x)$ we also deduce the analogous following symmetry for $\tilde H$:
$$\tilde H(x,p)= \tilde H(p,x)\quad \forall x,p\in U.$$
In the following we denote by $\frac{d}{dx}$ and $\frac{d}{dp}$ the (complex) derivative with respect to the first and the second variable of the function $\tilde H(\cdot,\cdot)$, respectively.
Therefore the Taylor expansion of $H_p(x)$ up to the order $m$ takes the form:
\beq\label{expa}H(x,p)-H(0,p)=H_p(x)-H_p(0)=
\sum_{k=1}^{m} \frac{1}{k!}\Re\Big(\frac{d^k\tilde H}{dx^k} (0,p) x^k\Big)+O(|x|^{m+1})\eeq uniformly for  $x\in\Omega$  and $p\in U$, 
where 
$$
\frac{d \tilde H }{dx}(0,p)=\frac{\partial^i  H_p }{dx_1}(0)-{\rm i} \frac{\partial^i  H_p }{dx_2}(0),$$
$$\frac{d^k \tilde H }{\partial x^k} (0,p) =\frac{d^k \tilde H_p }{dx^k}(0)=\frac{\partial^k  H_p }{dx_1^k}(0)-{\rm i} \frac{\partial^k  H_p }{dx_2dx_1^{k-1}}(0)\quad \forall k\geq 2.$$ 

\begin{lemma}\label{robin10} Using the same notation $b, \beta_i$ of the introduction, 
the following holds:
$$\begin{aligned}\sum _{i=0}^{\alpha-1} (H(x,\beta_i)-H(0,\beta_i))=&
\al \sum _{k=1}^{\al}\frac{1}{k!}\Re\Big(\frac{d^k\tilde H}{dx^k} (0,0)x^k\Big)+\frac{1}{(\al-1)!}\sum _{k=1}^{\al}\frac{1}{k!}\Re\Big(\frac{\partial^{k+\al}\tilde H}{\partial p^\al\partial x^k}(0,0)b\,x^k\Big)\\ &
+O(|b|^{2}|x|)+O(|x|^{\al+1}).\end{aligned}
$$
 uniformly for $b\in U$ and $x\in\Omega$. 
\end{lemma}
\begin{proof} According to \eqref{expa}  we compute
\beq\label{symmmm}\sum_{i=1}^\al\Big(H(x,\beta_i)-H(0,\beta_i)\Big)=\sum_{k=1}^{\al} \frac{1}{k!}\Re\Big(\sum_{i=1}^\al\frac{d^k\tilde H}{dx^k} (0,\beta_i) x^k\Big)+O(|x|^{\al+1}).\eeq
Let us expand the complex function $\frac{d^k\tilde H}{dx^k} (0,\beta_i)$:
$$\frac{d^k\tilde H}{dx^k} (0,\beta_i)=\frac{d^k\tilde H}{dx^k} (0,0)+\sum_{h=1}^{2\al-1}\frac{1}{h!}\frac{\partial^{k+h}\tilde H}{\partial p^h\partial x^k}(0,0)\beta_i^h+O(|\beta_i|^{2\al}).$$
Next we use that $ \sum _{i=0}^{\alpha-1} \beta_i^h=0$ for any $h$ which is not an integer multiple of $\alpha$, whereas $\beta_i^{\alpha j}=b^j$ for any integer $j$, by which
$$\sum_{i=0}^{\al-1}\frac{d^k\tilde H}{dx^k} (0,\beta_i)=\al \frac{d^k\tilde H}{dx^k} (0,0)+\frac{1}{(\al-1)!}\frac{\partial^{k+\al}\tilde H}{\partial p^\al\partial x^k}(0,0)b
+O(|b|^{2}). $$
By inserting the last identity into \eqref{symmmm} we obtain the thesis.

\end{proof}

\end{document}